\documentclass{gtpart}
\usepackage{amscd}
\usepackage[all]{xy}

\makeatletter\@addtoreset{equation}{section}
\makeatother
\makeatletter

\title{Gottlieb and Whitehead center groups of projective spaces}

\author{Marek Golasi\'nski}
\givenname{Marek}
\surname{Golasi\'nski}
\address{Faculty of Mathematics and Computer Science,\newline
Nicolaus Copernicus University,
87-100 Toru\'n,\newline Chopina 12/18, Poland}
\email{marek@mat.uni.torun.pl}
\author{Juno Mukai}
\givenname{Juno}
\surname{Mukai}
\address{Shinshu University,\
Matsumoto,\newline Nagano Pref.\ 390-8621, Japan}
\email{mukai@orchid.shinshu-u.ac.jp}



\newcommand{\rarrow}[1]{{\buildrel #1 \over \longrightarrow}}

\def\Ker{{\rm Ker }}
\def\Z{{\mathbb Z}}
\def\C{{\mathbb C}}
\def\R{{\mathbb R}}
\def\H{{\mathbb H}}
\def\K{{\mathbb K}}
\def\F{{\mathbb F}}

\def\S{{\mathbb S}}
\def\P{{\rm P}}
\def\dim{{\rm dim}}

\newtheorem{thm}{Theorem}[section]
\newtheorem{cor}[thm]{Corollary}
\newtheorem{prop}[thm]{Proposition}
\newtheorem{lem}[thm]{Lemma}
\newtheorem{rem}[thm]{Remark}

\newtheorem{exam}[thm]{Example}

\newtheorem{quest}[thm]{Question}
\newenvironment{pf}{\noindent{\bf Proof.}}

\begin{document}

\begin{abstract}
By use of Siegel's method and the classical results of homotopy groups of spheres and Lie groups,
we determine some Gottlieb groups of projective spaces or give the lower bounds of their orders.
Furthermore, making use of the properties of Whitehead products, we determine some  Whitehead center
groups of projective spaces.
\end{abstract}

\maketitle

\section*{Introduction}
Let $X$ be a connected and pointed space. The $k$-th {\em Gottlieb group}
$G_k(X)$ of $X$ has been defined in \cite{G1}. This is the subgroup of the $k$-th homotopy group
$\pi_k(X)$ containing all elements which can be represented by a map $f : \mathbb{S}^k\to X$
such that $f\vee\mbox{id}_X: \mathbb{S}^k\vee X\to X$ extends (up to homotopy) to a map
$F : \mathbb{S}^k\times X\to X$, where $\S^k$ is the $k$-sphere.
Throughout the paper we do not distinguish between a map and its homotopy class.
\par For $k\geq 1$, define the $k$-th {\em Whitehead center group} or {\em P-group} $P_k(X)\subseteq \pi_k(X)$ of elements $\alpha\in\pi_k(X)$
such that the Whitehead product $[\alpha,\beta]=0$ for all $\beta\in\pi_l(X)$ and $l\ge 1$.
Then, $G_k(X)\subseteq P_k(X)$ for all $k\geq 1$ \cite{G1}.
Furthermore, if $X$ is a Hopf space, $G_k(X)=P_k(X)=\pi_k(X)$ for all $k\ge 1$.

\par Let $\iota_n$ be the generator of $\pi_n(\mathbb{S}^n)$
represented by the identity map of $\mathbb{S}^n$. Given $\alpha\in\pi_k(\mathbb{S}^n)$,
$\alpha\in G_k(\mathbb{S}^n)$  if and only if $[\alpha,\iota_n]=0$.
We
note that $G_k(\mathbb{S}^n)=P_k(\mathbb{S}^n)$ \cite{GM}, \cite{V}.

\par Let $\mathbb{R}$ and $\mathbb{C}$ be the fields of reals and complex numbers, respectively and $\mathbb{H}$ the skew $\mathbb{R}$-algebra
of quaternions. Denote by $\mathbb{F}P^n$ for the $n$-projective space over $\mathbb{F}$.
Put $d=\dim_\mathbb{R}\mathbb{F}$, write
$i_{k,n,\mathbb{F}P}: \mathbb{F}P^k\hookrightarrow\mathbb{F}P^n$ for $k\le n$ the inclusion map
and write simply $i_\mathbb{F}=i_{1,n,\mathbb{F}P}: \mathbb{S}^d\hookrightarrow\mathbb{F}P^n$.

\par The purpose of this note is to determine some Gottlieb and $P$-groups of $\mathbb{F}P^n$. Our idea is
based on Lang's result $n!\pi_{2n+1}(\C P^n)\subseteq G_{2n+1}(\C P^n)$ \cite{L},
Barratt-James-Stein's result about the Whitehead products $[ -,\ i_\mathbb{F}]$
\cite{BJS}.

The method is to use Siegel's
result \cite{Si},
fibrations with $\mathbb{S}^{d(n+1)-1},\mathbb{F}P^n$ as base spaces and the results on the homotopy groups of spheres \cite{Mi0}, \cite{MT}, \cite{MMO}, \cite{T}, \cite{T2}, the classical
groups \cite{K}, \cite{ME}, \cite{Ma}, \cite{Ma1}, \cite{Mi}, \cite{Mi1}, \cite{MT1}, \cite{MT2}, \cite{Mo1}, \cite{Mo2}.

For example, $G_{k+n}(\mathbb{R}P^n)$ is partly determined and
$P_{k+n}(\mathbb{R}P^n)$  completely determined for $n\geq 2$ and $k\leq 7$.

We set
$$
P'_k(\mathbb{F}P^n)=P_k(\mathbb{F}P^n)\cap({\gamma_n}_\ast\pi_k(\S^{d(n+1)-1}))
$$
and
$$
P''_k(\mathbb{F}P^n)=P_k(\mathbb{F}P^n)\cap{i_\F}_\ast (E\pi_{k-1}(\S^{d-1})).
$$
Notice that $P''_k(\mathbb{F}P^n)=0$ if $d=1,2$ and $k\ge d+1$.
We determine the fact that
$P_{4n+3+k}(\H P^n)=P'_{4n+3+k}(\H P^n)\oplus P''_{4n+3+k}(\H P^n)$ and the group
$P'_{4n+3+k}(\H P^n)$ for $0\le k\le 10$.

Some particular cases of our results about the Gottlieb groups of $\mathbb{F}P^n$ overlap with those of \cite{LMW} and \cite{PW1}.

\section{Whitehead center groups of projective spaces}
Let $EX$ be the suspension of a space $X$ and denote by
$E: \pi_k(X)\to\pi_{k+1}(EX)$ the suspension
homomorphism. Write
$\eta_2\in\pi_3(\S^2)$, $\nu_4\in\pi_7(\S^4)$ and $\sigma_8
\in\pi_{15}(\S^8)$ for the appropriate Hopf maps, respectively. We set
$\eta_n=E^{n-2}\eta_2\in\pi_{n+1}(\S^n)$ for $n\ge 2$,
$\nu_n=E^{n-4}\nu_4\in\pi_{n+3}(\S^n)$ for $n\ge 4$ and
$\sigma_n=E^{n-8}\sigma_8\in\pi_{n+7}(\S^n)$ for $n\ge 8$. Note that $\nu_n$ for $n\geq 5$ and $\sigma_n$ for $n\geq 9$ stand for generators of $\pi_{n+3}(\S^n)\cong\Z_{24}$ and $\pi_{n+7}(\S^n)\cong\Z_{240}$, respectively. By abuse of notaion, these are also used to represent generators of the $2$-primary components $\pi^n_{n+3}\cong\Z_8$ and $\pi^n_{n+7}\cong\Z_{16}$, respectively.

Denote by $\gamma_n=\gamma_{n,\mathbb{F}} : \mathbb{S}^{(n+1)d-1}\to\mathbb{F}P^n$ the quotient map and by $q_n=q_{n,\mathbb{F}}
: \mathbb{F}P^n\to\mathbb{S}^{dn}$ the map pinching $\mathbb{F}P^{n-1}$ to the base point.
Then, as it is well known,

\begin{equation}\label{qgam}
q_n\gamma_n=\left\{ \begin{array}{ll} (1-(-1)^n)\iota_n,
& \mbox{if} \hspace{3mm} \mathbb{F}=\mathbb{R},\\
n\eta_{2n}, & \mbox{if} \hspace{3mm} \mathbb{F}=\mathbb{C},\\
\pm n\nu_{4n},&\mbox{if} \hspace{3mm} \mathbb{F}=\mathbb{H}.
\end{array}
\right.
\end{equation}

Let $\Delta=\Delta_{\mathbb{F}P}: \pi_k(\mathbb{F}P^n)\to\pi_{k-1}(\mathbb{S}^{d-1})$ be the connecting map. By \cite{BJS},
$$
\Delta({i_\mathbb{F}}_\ast E)=\mbox{id}_{\pi_{k-1}(\mathbb{S}^{d-1})}
$$
and
$$
\pi_k(\mathbb{F}P^n)={\gamma_n}_\ast\pi_k(\mathbb{S}^{d(n+1)-1})\oplus {i_\mathbb{F}}_\ast E\pi_{k-1}(\mathbb{S}^{d-1}).
$$

\par By \cite[Theorem (2.1)]{HW},
\begin{equation}\label{hw}
[\xi\circ\alpha,\eta\circ\beta]=0\ \mbox{if}\ [\xi,\eta]=0 \
(\xi\in\pi_m(X), \eta\in\pi_n(X), \alpha\in\pi_k(\S^m),
\beta\in\pi_l(\S^n)).
\end{equation}

According to \cite{V}, a map $f: X\to Y$ is {\it cyclic} if $f\vee \mbox{id}_Y: X\vee Y\to Y$ extends to $F: X\times Y\to Y$. The extension $F$ is called the
{\em associated map} with $f$ \cite{L}. We recall from \cite[Lemma 2.1]{Si} and \cite[Lemma 1.3]{V}:

\begin{lem}{\bf (Siegel-Varadarajan)} \label{SV}
Let $f: X\to Y$ be cyclic. Then, the composite $f\circ g: W\to Y$ is cyclic for any map $g: W\to X$.
\end{lem}

The following is obtained from (\ref{hw}) and Lemma \ref{SV}:
\begin{lem} \label{eg} $P_n(X)\circ\pi_k(\S^n)\subseteq P_k(X)$ \ and \ $G_n(X)\circ\pi_k(\S^n)\subseteq G_k(X)$.
\end{lem}

\par By (\ref{hw}), \cite[Corollary(7.4)]{BB} and \cite[(4.1-3)]{BJS},
we obtain a key formula determining the Whitehead center groups of $\mathbb{F}P^n$.

\begin{lem} \label{Y}
Let $h_0\alpha\in\pi_k(\mathbb{S}^{2d(n+1)-3})$ be the $0$-th Hopf-Hilton invariant for $\alpha\in\pi_k(\mathbb{S}^{d(n+1)-1})$.
Then:

\vspace{2mm}

$\mbox{\em (1)}\;\;[\gamma_n\alpha,i_\mathbb{R}]
=\left\{ \begin{array}{ll} 0,
& if\, n\, is\, \mbox{odd},\\
(-1)^k\gamma_n(-2\alpha+[\iota_n,\iota_n]\circ h_0\alpha), & if\, n\, is\, \mbox{even};
\end{array}
\right.
$

\vspace{2mm}

$\mbox{\em (2)}\;\;[\gamma_n\alpha,i_\mathbb{C}]=\left\{ \begin{array}{ll} 0,
&if\, n\,is\,\mbox{odd},\\
\gamma_n(\eta_{2n+1}\circ E\alpha+[\iota_{2n+1},\eta_{2n+1}]\circ Eh_0\alpha), &if\, n\,is\,\mbox{even};
\end{array}
\right.
$

\vspace{2mm}

$\mbox{\em (3)}\;\;[\gamma_n\alpha,i_\mathbb{H}]=\pm(n+1)\gamma_n(\nu_{4n+3}\circ E^3\alpha+[\iota_{4n+3},\nu_{4n+3}]\circ E^3h_0\alpha).$
\end{lem}

Hereafter, we use the results and notations of \cite{T} freely. Let $\nu'$ and $\alpha_1(3)$ be the generators of the $2$-primary $\pi^3_6$ and $3$-primary $\pi_6(\S^3;3)$
components of $\pi_6(\mathbb{S}^3)\cong\Z_{12}$, respectively.
Set $\omega=\nu'+\alpha_1(3)$ and $\alpha_1(n)=E^{n-3}\alpha_1(3)$ for $n\geq 3$.
We use the relations
\begin{equation}\label{Wio44}
[\iota_4,\iota_4]=2\nu_4-E\omega
\end{equation}
and $E^{n-3}\omega=2\nu_n$ for $n\geq 5$.

We write
$$
O_\mathbb{F}(n)=\left\{ \begin{array}{ll} O(n),
& \mbox{for} \hspace{3mm} \mathbb{F}=\mathbb{R};\\
U(n), & \mbox{for} \hspace{3mm} \mathbb{F}=\mathbb{C};\\
Sp(n),&\mbox{for} \hspace{3mm} \mathbb{F}=\mathbb{H},
\end{array}
\right.
$$
where $O(n)$, $U(n)$ and $Sp(n)$ denote the Lie groups of the orthogonal, unitary and symplectic $n\times n$ matrices,
respectively. Furthermore, we set
$$
SO_\mathbb{F}(n)=\left\{ \begin{array}{ll} SO(n),
& \mbox{for} \hspace{3mm} \mathbb{F}=\mathbb{R};\\
SU(n), & \mbox{for} \hspace{3mm} \mathbb{F}=\mathbb{C};\\
Sp(n),&\mbox{for} \hspace{3mm} \mathbb{F}=\mathbb{H}.
\end{array}
\right.
$$
Denote by $i_n=i_{n,\mathbb{F}}: SO_{\mathbb{F}}(n-1)\hookrightarrow SO_{\mathbb{F}}(n)$
and $p_n=p_{n,\mathbb{F}} : SO_{\mathbb{F}}(n)\to\mathbb{S}^{dn-1}$ the inclusion and projection maps, respectively.
Then, we consider the exact sequence induced from the fibration $SO_{\mathbb{F}}(n+1)\stackrel{SO_{\mathbb{F}}(n)}{\longrightarrow}\mathbb{S}^{d(n+1)-1}$:
$$
({\mathcal{SF}}^n_k) \ \cdots\rarrow{}
\pi_{k+1}(\mathbb{S}^{d(n+1)-1})\rarrow{\Delta_\mathbb{F}}\pi_k(SO_{\mathbb{F}}(n))\rarrow{i_\ast} \pi_k(SO_{\mathbb{F}}(n+1))\rarrow{p_\ast}\cdots,
$$
where $i=i_{n+1,\mathbb{F}}$ and $p=p_{n+1,\mathbb{F}}$.

Denote by $r: U(n)\to SO(2n)$, $c: Sp(n)\to SU(2n)$ the canonical map.
Let $J=J_\mathbb{R} : \pi_k(SO(n))\to \pi_{k+n}(\mathbb{S}^n)$ be the $J$-homomorphism,
$J_\mathbb{F} : \pi_k(SO_\mathbb{F}(n))\to \pi_{k+dn}(\mathbb{S}^{dn})$
the complex or symplectic $J$-homomorphism defined as follows: $J_\mathbb{C}=J\circ r_\ast: \pi_k(SU(n))\to\pi_{k+2n}(\mathbb{S}^{2n})$
and $J_\mathbb{H}=J\circ r_\ast\circ c_\ast: \pi_k(Sp(n))\to\pi_{k+4n}(\mathbb{S}^{4n})$.
Then, we see that
\begin{equation}\label{JDel}
E^{d-1}\circ J_\mathbb{F}\circ\Delta_\mathbb{F}=J\circ\Delta\;\; \mbox{for}\;\; \Delta=\Delta_\mathbb{R}.
\end{equation}
\par Let $\omega_{n,\mathbb{R}}\in\pi_{n-1}(O(n))$, $\omega_{n,\C}\in\pi_{2n}(U(n))$ and
$\omega_{n,\H}\in\pi_{4n+2}(Sp(n))$ be the characteristic elements for appropriate bundles.
We note that $\omega_{1,\H}=\omega$ is the Blakeres-Massey's element, $\omega_{n,\R}=\Delta(\iota_n)$,
$\omega_{n,\C}=\Delta_\C(\iota_{2n+1})$ and $\omega_{n,\H}=\Delta_\H(\iota_{4n+3})$.
\par As it is well-known,
\begin{equation}\label{char}
i_{2n+1,\R}r\omega_{n,\C}=\omega_{2n+1,\mathbb{R}}
\;\;\mbox{and}\;\;
i_{2n+1,\C}c\omega_{n,\H}=\omega_{2n+1,\C}.
\end{equation}
If $\Delta\alpha=0$ for $\alpha\in\pi_k(\S^n)$, then there exists a lift $[\alpha]\in\pi_k(SO(n+1))$ of $\alpha$.
For the inclusion
$$
i_{m,n}=i_{m,n,\R}=i_{n,\R}\circ\cdots\circ i_{m+1,\R}: SO(m)\hookrightarrow SO(n)\ (m\leq n-1),
$$
we set
$[\alpha]_n={i_{m,n}}_\ast[\alpha]\in\pi_k(SO(n))$, where  $[\alpha]\in\pi_k(SO(m))$ is a lift of $\alpha\in\pi_k(\S^{m-1})$.
Notice that $\pi_3(SO(4))=\{[\eta_2]_4,[\iota_3]\}\cong\Z^2$ and
any element of $\pi_k(SO(4))\cong\pi_k(SO(3))
\oplus\pi_k(\mathbb{S}^3)$
is uniquely represented by two elements of $\pi_k(\mathbb{S}^3)$:
$$
[\eta_2]_4\alpha+[\iota_3]\beta \hspace{.3cm} \mbox{for} \hspace{.3cm} \alpha,\beta\in\pi_k(\mathbb{S}^3).
$$
Notice that we can take $J[\iota_3]=\nu_4$ and $J[\eta_2]=\omega$.

Define the subgroups of $\pi_k(\S^{d(n+1)-1})$ and $\pi_{k-1}(\S^{d-1})$ as follows:
$$
M_k(\S^{d(n+1)-1})=M_{k,\F}(\S^{d(n+1)-1})=\{\alpha\in\pi_k(\S^{d(n+1)-1})\mid[\gamma_n\alpha,i_\F]=0\},
$$
$$
L'_{k-1}(\S^{d-1})=L'_{k-1,n}(\S^{d-1})=\{\beta\in\pi_{k-1}(\S^{d-1})\mid[i_{\F}E\beta,i_\F]=0\},
$$
$$
L''_{k-1}(\S^{d-1})=L''_{k-1,n}(\S^{d-1})=\{\beta\in\pi_{k-1}(\S^{d-1})\mid[i_{\F}E\beta,\gamma_n]=0\}
$$
and
$$
L_{k-1}(\S^{d-1})=L_{k-1,n}(\S^{d-1})=L'_{k-1,n}(\S^{d-1})\cap L''_{k-1,n}(\S^{d-1}).
$$

We also define the
subgroup $Q_k(\mathbb{S}^3)$ of $\pi_k(\mathbb{S}^3)$ by
$$
Q_k(\mathbb{S}^3)=\{\beta\in\pi_k(\mathbb{S}^3)\mid \langle\iota_3,\beta\rangle=0\},
$$
where $\langle -,\ - \rangle$ stands for the Samelson product.

For an abelian group $G$, we denote by $(G; p)$ its $p$-component. For example, $(\pi_k(\S^n);p)=\pi_k(\S^n;p)$.

We write
$$
'L^{d-1}_{k-1}=L'_{k-1}(\S^{d-1};2),\ ''L^{d-1}_{k-1}=L''_{k-1}(\S^{d-1};2),\
L^{d-1}_{k-1}=L_{k-1}(\S^{d-1};2),
$$
$$
P^{d(n+1)-1}_k=P_k(\S^{d(n+1)-1};2),\
Q^{d-1}_{k-1}=Q_{k-1}(\S^{d-1};2),\
M^{d(n+1)-1}_k=M_k(\S^{d(n+1)-1};2).
$$
Write $(-,-)$ for the greatest common divisor.
\begin{exam}\label{e2}{\em
$M_{4n+3,\H}(\S^{4n+3})=\frac{24}{(24,n+1)}\pi_{4n+3}(\S^{4n+3})$;\
$L''_{3,n}(\S^3)=\frac{24}{(24,n+1)}\pi_3(\S^3)$\ for\ $n\geq 1$;\
$Q_3(\S^3)=12\pi_3(\S^3)$.}
 \end{exam}

We show the following, in which (6) was suggested by K. Morisugi:
\begin{lem}\label{LQ}
\mbox{\em (1)}
$M_k(\S^{d(n+1)-1})\circ\pi_m(\S^k)\subseteq  M_m(\S^{d(n+1)-1})$,\\
$L'_{k-1}(\S^{d-1})\circ\pi_{m-1}(\S^{k-1})\subseteq L'_{m-1}(\S^{d-1})$, \\
$L''_{k-1}(\S^{d-1})\circ\pi_{m-1}(\S^{k-1})\subseteq L''_{m-1}(\S^{d-1})$,\\
$L_{k-1}(\S^{d-1})\circ\pi_{m-1}(\S^{k-1})\subseteq L_{m-1}(\S^{d-1})$ and $Q_k(\mathbb{S}^3)\circ\pi_m(\S^k)\subseteq Q_m(\mathbb{S}^3)$;

\mbox{\em (2)}
$L''_{k-1,n}(\S^3)=E^{-4n-3}(\Ker\ (n+1){\nu_{4n+3}}_\ast)$, where\\
${\nu_{4n+3}}_\ast: \pi_{k+4n+2}(\S^{4n+6})\to\pi_{k+4n+2}(\S^{4n+3})$ is the induced homomorphism. In particular,
$L'_{k-1,n}(\S^3)\subseteq L''_{k-1,n}(\S^3)$\
for \ $n$\ odd;

\mbox{\em(3)}
$Q_{k-1}(\S^3)=L'_{k-1,n}(\S^3)$ for $n\geq 2$;

\mbox{\em(4)}
Let $n\ge 2$. Then,
$L'_{k-1}(\S^3;3)\subseteq L''_{k-1}(\S^3;3)$ and\\ $L'_{k-1}(\S^3;p)=L''_{k-1}(\S^3;p)=\pi_{k-1}(\S^3;p)$ for an odd prime $p\geq 5$;

\mbox{\em(5)}
${i_\H}_\ast P_k(\S^4)\subseteq{i_\H}_\ast EQ_{k-1}(\S^3)$;\ Moreover,\\
${i_\H}_\ast P_k(\S^4)={i_\H}_\ast EQ_{k-1}(\S^3)$ provided that $2E^4Q_{k-1}(\S^3)=0$;

\mbox{\em(6)}
$L''_{k-1,n}(\S^3)=\pi_{k-1}(\S^3)$\ for $n\geq 1$\ if \ $5\leq k\leq 4n+2$.
\end{lem}
\begin{pf}
Since $\langle\iota_3,\beta\circ\delta\rangle=\langle\iota_3,\beta\rangle\circ E^3\delta$ \cite[(6.3)]{J4},
the last of (1) follows. The rest of (1) is a direct consequence of (\ref{hw}).

By the properties of Whitehead products \cite{WG} and Lemma \ref{Y}(3),
$$
(-1)^k[i_{\H}E\beta,\gamma_n]=[\gamma_n,i_{\H}E\beta]
=[\gamma_n,i_{\H}]\circ E^{4n+3}\beta
=\pm(n+1)\gamma_n\nu_{4n+3}E^{4n+3}\beta.
$$
This leads to the first of (2).
If $[E\beta,\iota_4]=0$, then $[E\beta,\nu_4]=0$ (\ref{hw}). This leads to the second of (2) for $n=1$.
Since $i_{\H}\circ\nu_4=0$ for $n\geq 2$ and $[E\beta,\iota_4]=(2\nu_4-E\omega)\circ E^4\beta$, the relation
$i_\mathbb{H}[E\beta,\iota_4]=0$ implies $i_\mathbb{H}E(\omega E^3\beta)=0$, and so $\omega E^3\beta=0$. Hence,
$$
(\triangle) \hspace{1.5cm} 2\nu_{4n+3}E^{4n+3}\beta=0
$$
and $(n+1)\nu_{4n+3}E^{4n+3}\beta=0$ for odd $n\geq 3$. This leads to the second of (2) for $n\geq 2$.

Since
$\langle\iota_3,\beta\rangle=\langle\iota_3,\iota_3\rangle\circ E^3\beta=\omega E^3\beta$, we obtain (3).

$(\triangle)$ for $\beta\in\pi_{k-1}(\S^3;3)$ implies
$\beta\in L''_{k-1}(\S^3;3)$ and the first half of
(4). The second half is a direct consequence of the definition of $L'_{k-1}(\S^{d-1})$ and $L''_{k-1}(\S^{d-1})$.

We recall
$$
\pi_k(\S^4)=E\pi_{k-1}(\S^3)\oplus{\nu_4}_\ast\pi_k(\S^7).
$$
Suppose that $\beta=E\beta_1+\nu_4\beta_2\in P_k(\S^4)$ for $\beta_1\in\pi_{k-1}(\S^3)$ and $\beta_2\in\pi_k(\S^7)$. Then, by the relation $i_\H\beta=i_\H E\beta_1$, we obtain
$$
0=i_\H[\iota_4,\beta]=i_\H[\iota_4,E\beta_1]=(i_\H E)(\omega E^3\beta_1).
$$
This implies the inclusion $P_k(\S^4)\subseteq EQ_{k-1}(\S^3)\oplus{\nu_4}_\ast \pi_k(\S^7)$ and leads to the first half of (5). On the other hand, for $\beta\in Q_{k-1}(\S^3)$, we obtain $[\iota_4,E\beta]=2\nu_4E^4\beta=0$ by the assumption. This means
$EQ_{k-1}(\S^3)\subseteq P_k(\S^4)$ and leads to the second half of (5).

We have $\nu_{4n+3}E^{4n+3}\beta=J([\iota_3]_{4n+3}\beta)$.
By the fact that $\pi_{k-1}(\S^3)$ is finite for $k\neq 4$ and
$$
\pi_{k-1}(SO(4n+3))\cong\left\{ \begin{array}{lll}
\Z,&\;\mbox{if}\;\; k\equiv 0\ (\bmod\ 4);\\
\Z_2,&\;\mbox{if}\;\; k\equiv 1,2\ (\bmod\ 8);\\
0,&\;\mbox{if}\;\; k\equiv 3,5,6,7\ (\bmod\ 8)
\end{array}
\right.
$$
if $k\leq 4n+2$ \cite{Bo}, we have (6) for $k\not\equiv 1,2\ (\bmod\ 8)$.

Assume that $\nu_{4n+3}E^{4n+3}\beta\neq 0$ for $k\equiv 1\ (\bmod\ 8)$.
Then, in view of \cite[Lemma 2]{K}, we have $\nu_{4n+3}E^{4n+3}\beta=J([\iota_3]_{4n+3}\beta)=J(\delta'\eta_{k-2})$,
where $\delta'\in\pi_{k-2}(SO(4n+3))$ is a generator. Since
$J: \pi_{k-1}(SO(4n+3))\to\pi_{k+4n+2}(\S^{4n+3})$ is a monomorphism, we obtain $\delta'\eta_{k-2}=[\iota_3]_{4n+3}\beta$. By \cite[Theorem 1.1]{DM},
$\delta'\eta_{k-2}$ has the $SO(6)$-of-origin. This is a contradiction and leads to (6) for $k\equiv\ 1\ (\bmod\ 8)$.
\par By the parallel argument, we have the assertion for $k\equiv 2\ (\bmod\ 8)$
This leads to (6) and completes the proof.
\end{pf}

{\bf Another proof of Lemma \ref{LQ}(6).}\
To consider $k\equiv 1,2\ (\bmod\ 8)$, write (according to Adams' notations \cite{A}) $j_r$ for the image of the generator in $\pi_r(SO)$ under
$J : \pi_r(SO)\to \pi_r^S$. First, assume that $\nu_{4n+3}E^{4n+3}\beta\neq 0$ for $k\equiv 1\ (\bmod\ 8)$.
Then, in view of \cite[Lemma 2]{K}, we have $\nu_{4n+3}E^{4n+3}\beta=\nu E^\infty\beta=j_{k-2}\eta$
and $j_{k-2}\eta^2$ gene\-rates the $J$-image of $\pi_k(SO)\cong\pi_k(SO(4n+3))\cong\Z_2$. On the other hand, by
\cite[Proposition 3.1]{T}, $\beta\wedge\eta_2=E^2(\beta\eta_{k-1})=\eta_5E^3\beta$. This and the relation $\nu\eta=0$
imply a contradictory relation
$j_{k-2}\eta^2=\nu E^\infty(\beta\eta)=\nu\eta E^\infty\beta=0$.

Now, assume that $\nu_{4n+3}E^{4n+3}\beta\neq 0$ for $5\leq k\leq 4n$ with $k\equiv 2\ (\bmod\ 8)$.
Then, by \cite[Lemma 2]{K}, $\nu_{4n+3}E^{4n+3}\beta=\nu E^\infty\beta=\nu E^\infty\beta=j_{k-1}$.
By the relation (\cite[Proposition 3.1]{T})
\begin{equation}\label{nuwdb}
\nu_4\wedge\beta=\nu_7E^7\beta=(E^4\beta)\nu_{k+3},
\end{equation}
we get for the stable Toda bracket
$$
<j_{k-1}, \eta, 2\iota>=<\nu E^\infty\beta,\eta,2\iota>=<(E^\infty\beta)\nu,\eta,2\iota>$$
$$
\supseteq(E^\infty\beta)<\nu,\eta,2\iota>\equiv\, 0\ (\bmod\ 2\pi^S_{k-1}).
$$
On the other hand, by \cite[Corollary 11.7]{A} and \cite[Theorem B.\ v)]{M}, we deduce that $\pi^S_{k-1}$
contains the direct summand $\Z_8$ which is generated by $<j_{k-1},j_3, 24\iota>.$
(We point out that \cite[Theorem B.\ v)]{M} has a misprint.
It should be  $e(\rho_j)=2^{-p(j)}\ (\bmod\; 2^{-p(j)-1})$.)

Because, $j_1= \eta$, $12j_3= \eta^3$ and $j_{k-1}=j_{k-3}\circ\eta^2$, we derive
$$<j_{k-3}, j_3, 24\iota> \subseteq<j_{k-3},12j_3,2\iota>
=<j_{k-3}, \eta^3, 2\iota> \supseteq <j_{k-1},\eta,2\iota>.$$
Thus,
$$<j_{k-3}, j_3, 24\iota>\equiv <j_{k-1},\eta,2\iota> (\bmod\,2 \pi^S_{8s+3}).$$
This is a contradiction with the above which leads to the assertion.
\bigskip

Now, we show:
\begin{prop}\label{newP}
$P'_k(\mathbb{F}P^n)={\gamma_n}_\ast(P_k(\mathbb{S}^{d(n+1)-1})\cap M_k(\mathbb{S}^{d(n+1)-1}))$ and\\  $P''_k(\mathbb{F}P^n)={i_\F}_\ast (EL_{k-1}(\S^{d-1}))$.
\end{prop}
\begin{pf}
For an element $\alpha\in\pi_k(\S^{d(n+1)-1})$, it holds that $[\gamma_n\alpha,\gamma_n]=0$ if and only if $[\alpha,\iota_{d(n+1)-1}]=0$. This leads to the first half.

Next, for an element $\beta\in\pi_{k-1}(\S^{d-1})$, it holds that ${i_\F}E\beta\in P''(\F P^n)$ if and only if
$\beta\in L_{k-1}(\S^{d-1})$. This leads to the second half and completes the proof.
\end{pf}

We recall the order of $[\iota_n,\alpha]$ for $\alpha=\iota_n, \eta_n,\eta^2_n,E\omega\ (n=4),\nu_n, \nu^2_n, E\sigma'\ (n=8)$ and $\sigma_n$ \cite{GM}.

\begin{equation} \label{W1}
\sharp[\iota_n,\iota_n]=\left\{ \begin{array}{ll}
1, & \mbox{if}\;\;n=1,3,7,\\
2, & \mbox{if}\;\;n\;\;\mbox{is odd}\;\;\mbox{and}\;\;n\not=1,3,7,\\
\infty, & \mbox{if}\;\;n\;\;\mbox{is even};
\end{array}
\right.
\end{equation}

\begin{equation} \label{W2}
\sharp[\iota_n,\eta_n]=\left\{ \begin{array}{ll} 1,
& \mbox{if}\;\; n\equiv 3\ (\bmod\,4),\; n=2 \
\mbox{or}\ n=6,\\
2, & \mbox{if otherwise};
\end{array}
\right.
\end{equation}

\begin{equation} \label{W3}
\sharp[\iota_n,\eta^2_n]=\left\{ \begin{array}{ll} 1,
& \mbox{if}\;\;n\equiv 2,3\ (\bmod\,4),\\
2, & \mbox{if otherwise};
\end{array}
\right.
\end{equation}

\begin{equation} \label{W40}
\sharp[\iota_4,E\omega]=6;
\end{equation}

\begin{equation} \label{W4}
\sharp[\iota_n,\nu_n]=\left\{\begin{array}{ll}
1,&\quad\mbox{if}\;n\equiv 7\ (\bmod\,8),\ n= 2^i - 3\geq 5,\\
2,&\quad\mbox{if} \; n\equiv 1,3,5\ (\bmod\,8)\geq 9,\ n\neq 2^i - 3,\\
12,&\quad\mbox{if} \;n\equiv 2\ (\bmod\,4)\geq 6, n=4, 12,\\
24,&\quad \mbox{if} \;n\equiv 0\ (\bmod\,4)\geq 8,\ n\neq 12;
\end{array}
\right.
\end{equation}

\begin{equation} \label{W5}
\sharp[\iota_n,\nu_n^2]=1 \ \mbox{if and only if} \ n\equiv 4,5,7\ (\bmod\,8)\ \mbox{or} \ n = 2^i - 5\ \mbox{for} \ i\geq 4;
\end{equation}

\begin{equation}\label{W60}
\sharp[\iota_8,E\sigma']=60;
\end{equation}

\begin{equation}\label{W6}
\sharp[\iota_n,\sigma_n]=\left\{\begin{array}{ll}
1,&\ \mbox{if} \ n=11,\ n\equiv 15\ (\bmod\,16),\\
2,&\ \mbox{if}\ n\equiv 1\ (\bmod\,2)\geq 9,\ n\neq 11,\ n\not\equiv 15\ (\bmod\,16),\\
120,&\ \mbox{if} \ n=8,\\
240,&\ \mbox{if} \ n\equiv 0\ (\bmod\,2)\geq 10.
\end{array}
\right.
\end{equation}

We show:
\begin{lem}\label{PM}
\mbox{\em (1)}
$M_k(\S^{d(2n+2)-1})=\pi_k(\S^{d(2n+2)-1})$ for $\F=\R,\C$;

\mbox{\em (2)}\
\mbox{\em (i)}
$E\pi_{k-1}(\S^{2n-1})\cap\Ker\ {2\iota_{2n}}_\ast\subseteq M_{k,\R}(\S^{2n})$, where
${2\iota_{2n}}_\ast: \pi_k(\S^{2n})\to\pi_k(\S^{2n})$ is the induced homomorphism; \ $M_{k,\R}(\S^{2n})=\Ker\ {2\iota_{2n}}_\ast$ if $k\leq 4n-2$;

\hspace{.5cm}
\mbox{\em (ii)}
$E\pi_{k-1}(\S^{4n})\cap E^{-1}(\Ker\ {\eta_{4n+1}}_\ast)\subseteq M_{k,\C}(\S^{4n+1})$, where
${\eta_{4n+1}}_\ast: \pi_{k+1}(\S^{4n+2})\to\pi_{k+1}(\S^{4n+1})$
; \
$M_{k,\C}(\S^{4n+1})=E^{-1}(\Ker\ {\eta_{4n+1}}_\ast)$ \ if $k\leq 8n$;

\hspace{.5cm}
\mbox{\em (iii)}
$E\pi_{k-1}(\S^{4n+2})\cap E^{-3}(\Ker\ (n+1){\nu_{4n+3}}_\ast)\subseteq M_{k,\H}(\S^{4n+3})$, where\
${\nu_{4n+3}}_\ast: \pi_{k+3}(\S^{4n+6})\to\pi_{k+3}(\S^{4n+3})$; \
$M_{k,\H}(\S^{4n+3})=E^{-3}(\Ker\ (n+1){\nu_{4n+3}}_\ast)$ \ if $k\leq 8n+4$.
In particular,
$M_{k,\H}(\S^{4n+3})=\pi_k(\S^{4n+3})$ if $n+1\equiv 0\ (\bmod\ 24)$;

\hspace{.5cm}
\mbox{\em (iv)}
$E\pi_{k-1}(\S^{8n+2})\cap E^{-3}(\Ker\ (2n+1){\nu_{8n+3}}_\ast)\subseteq M_{k,\H}(\S^{8n+3})$; \\
$M_{k,\H}(\S^{8n+3})=E^{-3}(\Ker\ (2n+1){\nu_{8n+3}}_\ast)$ \ if \ $k\leq 16n+4$;

\hspace{.5cm}
\mbox{\em (v)}
$M_{k,\H}(\S^{8n+7})=E^{-3}(\Ker\ 2(n+1){\nu_{8n+7}}_\ast)$;

\mbox{\em (3)}
$[\iota_{2n},\iota_{2n}]\in M_{4n-1,\R}(\S^{2n})$;

\mbox{\em (4)}
$M_{k,\R}(\S^2)=\pi_k(\S^2)$ except $k=2$ and $M_{2,\R}(\S^2)=0$;

\mbox{\em (5)}
$P_{k+2n}(\S^{2n})\subseteq M_{k+2n,\R}(\S^{2n})$ \ if  $k\leq 2n-2$;

\mbox{\em (6)}
\mbox{\em (i)}
$P_k(\S^{8n+3};p)=M_{k,\H}(\S^{8n+3};p)
=\pi_k(\S^{8n+3};p)$ \ for an odd prime $p\geq 5$;

\hspace{.5cm}
\mbox{\em (ii)}
$M_{k,\H}(\S^{8n+3};3)
=\pi_k(\S^{8n+3};3)$ \ if $n\equiv 1\ (\bmod\ 3)$.

\mbox{\em (7)}
$M^{8n+3}_{k+8n+3}=\left\{\begin{array}{ll}
8\pi^{8n+3}_{8n+3},& if\, k=0;\\
\pi^{8n+3}_{k+8n+3},& if\;k=1,2,4,5,7,8,9,10;\\
2\pi^{8n+3}_{k+8n+3},& if\;k=3,6.
\end{array}
\right.$
\end{lem}
\begin{pf}
(1) is a direct consequence of Lemma \ref{Y}(1);(2).

By Lemma \ref{Y}(1), $E\alpha\in M_{k,\R}(\S^{2n})$ if and only if $2E\alpha=2\iota_{2n}\circ E\alpha=0$. This leads to the first half of (2)-(i). The second half is obtained from the Freudenthal suspension theorem.
By Lemma \ref{Y}(2), $E\alpha\in M_{k,\C}(\S^{4n+1})$ if and only if $\eta_{4n+1}E^2\alpha=0$ for $\alpha\in\pi_{k-1}(\S^{4n})$. This leads to (2)-(ii).

By Lemma \ref{Y}(3), $E\alpha\in M_{k,\H}(\S^{4n+2})$ if and only if $(n+1)\nu_{4n+3}E^4\alpha=0$. This leads to (2)-(iii). By the parallel argument, we have
(2)-(iv). By Lemma \ref{Y}(3) and the relation $[\iota_{8n+7},\nu_{8n+7}]=0$ (\ref{W4}), we obtain
(2)-(v).

(3) is a direct consequence of \cite[Proposition 2]{Hi1}.

By Lemma \ref{Y}(1), $[\gamma_2,i_{\R}]=-2\gamma_2$ and $[\gamma_2\eta_2,i_{\R}]=0$. This implies
$M_{2,\R}(\S^2)=0$ and
$\eta_2\in M_{3,\R}(\S^2)$. By Lemma \ref{LQ}(1),
$\pi_k(\S^2)=\eta_2\circ\pi_k(\S^3)\subseteq M_{k,\R}(\S^2)$. This leads to (4).

Suppose that $E\alpha\in\ P_{2n+k}(\S^{2n})$. Then, $0=[\iota_{2n},E\alpha]=[\iota_{2n},\iota_{2n}]\circ E^{2n}\alpha.$
By the fact that $H[\iota_{2n},\iota_{2n}]=\pm 2\iota_{4n-1}$ for the Hopf homomorphism
$H: \pi_{4n-1}(\S^{2n})\to\pi_{4n-1}(\S^{4n-1})$, we have $0=H[\iota_{2n},\iota_{2n}]\circ E^{2n}\alpha=2E^{2n}\alpha$.
Since
$E^{2n-1}: \pi_{k+2n}(\S^{2n})\to\pi_{k+4n-1}(\S^{4n-1})$
is an isomorphism, we obtain $2E\alpha=0$.
Hence, (2)-(i) leads to (5).

Since $s[\iota_{8n+3},\alpha]=0$ for $s=2,p$ and $t\alpha_1(8n+3)E^3\alpha=0$ for $t=3,p$. This leads to (6)-(i). The assumption on $n$ implies $2n+1\equiv 0\ (\bmod \ 3)$. We have $u[\iota_{8n+3},\nu_{8n+3}]\circ E^3h_0(\alpha)=0$ for $u=2,3$. This leads to (6)-(ii).

Next, we show (7). By Lemma \ref{Y}(3), we get that $\alpha\in M^{8n+3}_{k+8n+3}$ if and only
if $\nu_{8n+3}E^3\alpha=0$. Hence, $M^{8n+3}_{8n+3}=8\pi^{8n+3}_{8n+3}$ and in view of the relation $\nu_6\eta_9=0$ \cite[(5.9)]{T}, we get that
$M^{8n+3}_{k+8n+3}=\pi_{k+8n+3}^{8n+3}$ for $k=1,2$.
\par Because $\pi_{8n+6}^{8n+3}=\{\nu_{8n+3}\}$ and
$\sharp\nu^2_{8n+3}=2$,
we derive that $M^{8n+3}_{8n+6}=2\pi^{8n+3}_{8n+6}$.
\par 
Trivially, $M_{8n+7}^{8n+3}=M_{8n+8}^{8n+3}=0$.
\par 
Since $\sharp\nu^3_{8n+3}=2$ and $\pi_{8n+9}^{8n+3}=\{\nu^2_{8n+3}\}$,
we have $M_{8n+9}^{8n+3}=2\pi_{8n+9}^{8n+3}=0$.
\par Because $\pi_{8n+10}^{8n+3}=\{\sigma_{8n+3}\}$, in view of the relation $\nu_{11}\sigma_{14}=0$ \cite[(7.20)]{T},
we deduce $M_{8n+10}^{8n+3}=\pi_{8n+10}^{8n+3}$.
\par In virtue of the relation \cite[(7.17) and (7.18)]{T}
\begin{equation}\label{nu6ep}
\nu_6\varepsilon_9=\nu_6\bar{\nu}_9=[\iota_6,\nu^2_6],
\end{equation}
we obtain
$\nu_7\varepsilon_{10}=\nu_7\bar{\nu}_{10}=0$. So we get $M_{8n+11}^{8n+3}=\pi_{8n+11}^{8n+3}$.
Trivially, $M_{8n+12}^{8n+3}=\pi_{8n+12}^{8n+3}$.
\par Finally, \cite[(5.9)]{T} leads to $M^{8n+3}_{8n+13}=\pi^{8n+3}_{8n+13}$ and the proof is completed.
\end{pf}

By Lemma \ref{PM} and the fact that $\eta^3_5=12\nu_5$, $\eta_5\nu_6=0$, $\pi_{n+6}(\S^n)=\{\nu^2_n\}\cong\Z_2$ for $n\geq 5$,
we obtain:
\begin{exam}\label{ex1}
{\em $M_{7,\C}(\S^5)=0$,\ $M_{12,\C}(\S^9)=\pi_{12}(\S^9)$ \ and \
$M_{4n+6,\H}(\S^{4n+3})=\frac{2}{(2,n+1)}\pi_{4n+6}(\S^{4n+3})$.}
\end{exam}

Here, we show the following result in the abelian groups.

\begin{lem}\label{Mrk}
Let $G$ be a group with  $G=G_1\oplus G_2$ and write
$p_k : G\to G_k$ for the projections with $k=1,2$.
If $H$ is a subgroup of $G$ such that\\
$$
(\flat) \quad p_k(H)<H\ \mbox{for}\ k=1,2
$$
then $H=G_1\cap H\oplus G_2\cap H$.
\end{lem}
\begin{pf}
Because $p_1(H)<G_1$ and $p_2(H)<G_2$, we have
$H=p_1(H)\oplus p_2(H)$.
Hence, by the assumotion $p_k(H)<H$ for $k=1,2$, $p_1(H)<H\cap G_1$ and $p_2(H)<H\cap G_2$.
Certainly, $H\cap G_1<p_1(H)$ and $H\cap G_2<p_2(H)$.
Consequently, $H=H\cap G_1\oplus H\cap G_2$.
\end{pf}

As it is well known,
\begin{equation}\label{JDel2}
J(\Delta\alpha)=\pm[\iota_n,\alpha], \ \mbox{where} \  \alpha\in\pi_k(\mathbb{S}^n).
\end{equation}
For a subgroup $M\subseteq\pi_k(\mathbb{S}^n)$, we set $[\iota_n,M]=J(\Delta M)$.
Let $P: \pi_{k+2}(\S^{2n+1})\to\pi_k(\S^n)$ be the $P$ homomorphism in the EHP sequence defined as the notation ``$\Delta$" \cite{T}. Then, by \cite[Corollary (7.4)]{BB} and \cite[Proposition 2.5]{T},
$$
[\iota_n,\pi^n_k]=P(E^{n+1}\pi^n_k).
$$

Consider the homomorphisms $[- , i_\mathbb{F}]:\pi_k(\mathbb{F}P^n)\to\pi_{k+d-1}(\mathbb{F}P^n)$ and $[- , \gamma_n]:\pi_k(\mathbb{F}P^n)\to\pi_{k+d(n+1)-2}(\mathbb{F}P^n)$. By the additivity of the Whitehead products \cite{WG} and (\ref{hw}), we obtain
$$
P_k(\mathbb{F}P^n)=\Ker\ [ - , i_\mathbb{F}]\cap\Ker\ [ - , \gamma_n].
$$

Then, we show:
\begin{lem}\label{ker}
Let $n\geq 2$. Then:

\mbox{\em (1)}
$\Ker\ [ - , i_\mathbb{F}]={\gamma_n}_\ast M_k(\S^{d(n+1)-1})\oplus {i_{\F}}_\ast EL'_{k-1}(\S^{d-1})$;

\mbox{\em (2)}
For $\F=\R, \C$ and $k\geq d+1$:

\hspace{.2cm}
{\em (i)}
$\Ker\ [ - , i_\mathbb{F}]={\gamma_n}_\ast M_k(\S^{d(n+1)-1})$;

\hspace{.2cm}
{\em (ii)}
$\Ker\ [ - , \gamma_n]={\gamma_n}_\ast P_k(\S^{d(n+1)-1})$;

\mbox{\em (3)} For $\F=\H$:

\hspace{.2cm} {\em (i)}
$(\Ker\ [ - , \gamma_n] ;p)={\gamma_n}_\ast \pi_k(\S^{4n+3};p)\oplus{i_{\H}}_\ast EL''_{k-1}(\S^3;p)$ for an odd prime $p$;

\hspace{.2cm} {\em (ii)}
$(\Ker\ [ - , \gamma_n]; 2)={\gamma_n}_\ast P_k(\S^{4n+3};2)\oplus{i_{\H}}_\ast EL''_{k-1}(\S^3;2)$
provided
that \\
$$
\quad\quad (\ast) \quad\quad [\iota_{4n+3},\pi^{4n+3}_k]\cap(n+1)\nu_{4n+3}E^{4n+3}\pi^3_{k-1}=0.
$$
\end{lem}
\begin{pf}
(1) follows from Lemma \ref{Y} and the relation
$$
0=
[\gamma_n\alpha+i_\mathbb{F}E\beta,i_\mathbb{F}]
=[\gamma_n\alpha,i_\mathbb{F}]
+[i_\mathbb{F}E\beta,i_\mathbb{F}].
$$
In the quaternionic case, this is equivalent
the following two formulas:
\begin{equation}\label{iH}
(n+1)(\nu_{4n+3}E^3\alpha+[\iota_{4n+3},\nu_{4n+3}]\circ E^3h_0(\alpha))=0
\end{equation}
and
\begin{equation}\label{iH1}
\nu'E^3\beta=0.
\end{equation}

By the assumption of (2), $\pi_{k-1}(\S^{d-1})=0$. (2)-(i) follows from (1). (2)-(ii) is obtained by the fact that
$[\alpha,\iota_{d(n+1)-1}]=0$ if and only if $[\gamma_n\alpha,\gamma_n]=0$.

For an element $\gamma_n\alpha+i_\mathbb{H}E\beta\in\Ker\ [ - , \gamma_n]$, we see that
$$
0=[\gamma_n\alpha+i_\mathbb{H}E\beta,\gamma_n]
=\gamma_n[\alpha,\iota_{4n+3}]+
[i_\mathbb{H}E\beta,\gamma_n]
$$
$$
=\gamma_n([\alpha,\iota_{4n+3}]\pm(n+1)\nu_{4n+3}E^{4n+3}\beta).
$$
That is,
\begin{equation}\label{gmH}
[\alpha,\iota_{4n+3}]=\pm(n+1)\nu_{4n+3}E^{4n+3}\beta.
\end{equation}

Since $[\iota_{4n+3},\pi_k(\S^{4n+3};p)]=0$ and $P_k(\S^{4n+3};p)=\pi_k(\S^{4n+3};p)$ for an odd prime $p$,
(3)-(i) follows. The assumption $(\ast)$ in (3)-(ii) implies $\alpha\in P_k(\S^{4n+3};2)$ and $\beta\in L''_{k-1}(\S^3;2)$.
This leads to (3)-(ii).
\end{pf}

Notice that
the condition $(\ast)$ in (3)-(ii) implies
the assumption $(\flat)$ of Lemma \ref{Mrk}.

Now, we obtain:

\begin{prop}\label{rc1}
\mbox{\em (1)}
Let $\F=\R, \C$  and $k\ge d+1$. Then,
$P_k(\mathbb{F}P^n)=
{\gamma_n}_\ast(P_k(\mathbb{S}^{d(n+1)-1})\cap M_k(\S^{d(n+1)-1}))$\ for $n\geq 2$. In particular, $P_k(\mathbb{F}P^{2n+1})={\gamma_{2n+1}}_\ast P_k(\mathbb{S}^{d(2n+2)-1})$\ for $n\geq 1$;

\mbox{\em (2)}
$P_2(\mathbb{R}P^2)=0$, $P_k(\mathbb{R}P^2)=\pi_k(\mathbb{R}P^2)$ for $k\geq 3$ and\\
$P_{k+2n}(\mathbb{R}P^{2n})={\gamma_{2n}}_\ast P_{k+2n}(\mathbb{S}^{2n})$\ if\ $k\leq 2n-2$;

\mbox{\em (3)}
$P_{k+4n+1}(\mathbb{C}P^{2n})={{\gamma_{2n}}_\ast (P_{k+4n+1}(\mathbb{S}^{4n+1})\cap E^{-1}(\Ker\ {{\eta_{4n+1}}_\ast}}))$, if\\
$k\leq 4n-1$.
\end{prop}
\begin{pf}
By Lemmas \ref{ker}(2) and \ref{PM}(1), we obtain (1).

By (1), Lemma \ref{PM}(1);(4);(5) and the fact that $P_k(\S^2)=\pi_k(\S^2)$ for $k\geq 3$, we obtain (2). By (1) and Lemma \ref{PM}(2)-(ii), we obtain (3). This completes the proof.
\end{pf}

In particular, we derive that $P_k(\mathbb{C}P^3)=\pi_k(\mathbb{C}P^3)$ for all $k\ge 1$ \cite{L}.

In the quaternionic case, we obtain:

\begin{prop}\label{rc12}
\mbox{\em (1)}
$P_k(\mathbb{H}P^{2n+1})=
{\gamma_{2n+1}}_\ast(P_k(\mathbb{S}^{8n+7})\cap M_k(\S^{8n+7}))$\\
$\oplus{i_\H}_*EQ_{k-1}(\S^3)$. In particular, $P_k(\mathbb{H}P^{2n+1})={\gamma_{2n+1}}_\ast P_k(\mathbb{S}^{8n+7})\oplus{i_\mathbb{H}}_*EQ_{k-1}(\mathbb{S}^3)$, if $n+1\equiv 0\ (\bmod\ 12)$;

\mbox{\em (2)}
\mbox{\em (i)}
$P_k(\mathbb{H}P^{2n};3)=
{\gamma_{2n}}_\ast M_k(\S^{8n+3};3)\oplus{i_\H}_*EQ_{k-1}(\S^3;3)$.
In particular, $P_k(\mathbb{H}P^{2n};3)=
{\gamma_{2n}}_\ast\pi_k(\mathbb{S}^{8n+3};3)\oplus{i_\H}_*EQ_{k-1}(\S^3;3)$ if $n\equiv 1\ (\bmod\ 3)$;

\hspace{.4cm}
\mbox{\em (ii)}
$P_k(\mathbb{H}P^{2n};p)=
{\gamma_{2n}}_\ast\pi_k(\mathbb{S}^{8n+3};p)\oplus{i_\H}_*E\pi_{k-1}(\S^3;p)$ for an odd prime $p\geq 5$;

\hspace{.4cm}
\mbox{\em (iii)}
$P_k(\mathbb{H}P^{2n};2)={\gamma_n}_\ast (P^{8n+3}_k\cap M^{8n+3}_k)\oplus{i_\H}_*EL^3_{k-1}$ provided that
$$
\quad (\heartsuit) \quad [\iota_{8n+3},M^{8n+3}_k]\cap \nu_{8n+3}\circ E^{8n+3}{}'L^3_{k-1}=0;
$$

\mbox{\em (3)}
$P_k(\mathbb{H}P^n)=
{i_\H}_*EQ_{k-1}(\S^3)$ for $5\leq k\leq 4n+2$ 
and $n\ge 2$.
\end{prop}
\begin{pf}
(1) is a direct consequence of Lemmas \ref{LQ}(2),
\ref{ker}(1);(3), \ref{PM}(2)-(iii) and the fact that
$$
\Ker\ ([ - , \gamma_n]\mid_{\Ker\ [ - , i_\mathbb{F}]})=\Ker\ [ - , i_\mathbb{F}]\cap\Ker\ [ - , \gamma_n].
$$

Let $n$ be even. Then, Lemmas \ref{LQ}(4) and  \ref{ker}(3)-(i) lead to the first half of (2)-(i). The second half follows from Lemma \ref{PM}(6)-(ii).
Lemmas \ref{LQ}(4), \ref{PM}(6)-(i) and \ref{ker}(4)-(i) imply (2)-(ii). By Lemma \ref{ker}(1);(3)-(ii), we have (2)-(iii). By Lemmas \ref{ker}(1);(3) and \ref{LQ}(3);(6), we obtain (3). This completes the proof.
\end{pf}

\begin{cor}\label{alphal}
Let $n$ be even. Suppose that $\alpha\in\pi^{4n+3}_k$, $\beta\in\pi^3_{k-1}$ satisfy {\em (\ref{iH}), (\ref{iH1}) and (\ref{gmH})}.

\mbox{\em (1)}
Let $E\alpha'\in M^{4n+3}_l$, that is,
$\alpha'\in\pi^{4n+2}_{l-1}$ be an element satisfying $(n+1)\nu_{4n+3}E^4\alpha'=0$.
If $[\alpha,E\alpha']\ne 0$, then,
$\gamma_n\alpha+i_{\H}E\beta\not\in P_k(\H P^n)$;

\mbox{\em (2)}
Suppose that $\pi^{4n+3}_k$ is a cyclic group generated by $\alpha$ with $[\alpha,E\alpha']\ne 0$. Then, the condition $(\heartsuit)$ in Proposition {\em \ref{rc12}(2)-(iii)} holds.
\end{cor}
\begin{pf}
By (\ref{iH}), (\ref{gmH}) and (\ref{nuwdb}),
$$
[\alpha,E\alpha']=[\alpha,\iota_{4n+3}]\circ E^k\alpha'=
\pm(n+1)(E^{4n}\beta)\circ E^{k-4}(\nu_{4n+3}E^4\alpha')=0.
$$
This completes the proof.
\end{pf}

\section{Whitehead center groups of real and complex projective spaces}
In this section, we determine some $P$-groups of real and complex projective spaces.
By Lemma \ref{Y} and the fact that $\sharp[i_\F,i_\F]=1$ for $\F=\R,\C$, we have that
$P_1(\R P^n)=\frac{1+(-1)^{n-1}}{2}\pi_1(\R P^n)$ and
$P_2(\C P^n)=\frac{3+(-1)^n}{2}\pi_2(\C P^n)$.

Now, we show:
\begin{thm}\label{main0}
The equality $P_{k+n}(\mathbb{R}P^n)={\gamma_n}_\ast P_{n+k}(\mathbb{S}^n)$ holds if $k\le 7$ except
the pairs: $(k,n)=(4,4),(5,4),(6,4)$.
\par Furthermore, $P_8(\mathbb{R}P^4)=\{\gamma_4(E\nu')\eta_7\}\cong\Z_2$,\
$P_9(\mathbb{R}P^4)=\{\gamma_4(E\nu')\eta^2_7\}\cong\Z_2$ and $P_{10}(\mathbb{R}P^4)
=\{3\gamma_4\nu^2_4, \gamma_4(\nu_4-\alpha_1(4))\alpha_1(7)\}\cong\Z_{24}$.
\end{thm}
\begin{pf}
For odd $n$, the equality holds by Proposition \ref{rc1}(1). For $n=2$ or the stable case, the equality
holds by Proposition \ref{rc1}(2). So, it suffices to prove
the
pairs $(k,n)=(3,4), (4,4), (5,4), (5,6), (6,4), (6,6), (7,4), (7,6)$ and $(7,8)$.

First, we calculate $P_7(\mathbb{R}P^4)$.
By Lemma \ref{Y}(1) and (\ref{Wio44}), we obtain
\begin{equation}\label{iRnu}
[\gamma_4\nu_4,i_\R]=\gamma_4E\omega.
\end{equation}
By Lemma \ref{PM}(3), $[\iota_4,\iota_4]\in M_7(\S^4)$. In virtue of
Lemma \ref{Y}(1), $[\gamma_4E\omega,i_\R]=-2\gamma_4E\omega$. So,
$M_7(\S^4)=\{12\nu_4,[\iota_4,\iota_4], 6E\omega\}=P_7(\mathbb{S}^4)$ \cite{GM} and hence, Proposition \ref{rc1}(1) leads to the equality
$P_7(\mathbb{R}P^4)={\gamma_4}_\ast P_7(\mathbb{S}^4)$.
By the parallel argument, we obtain $P_{15}(\mathbb{R}P^8)={\gamma_8}_\ast P_{15}(\mathbb{S}^8)$.

Next, we determine $P_k(\mathbb{R}P^4)$ for $k=8,9,10$. We recall $P_k(\S^4)=\pi_k(\S^4)$ for $8\leq k\leq 10$ \cite{GM}. By Lemma \ref{PM}(2)-(i) and (\ref{iRnu}), $(E\nu')\eta_7\in M_8(\S^4)$ and
$$
[\gamma_4\nu_4\eta_7,i_\R]=\gamma_4(E\nu')\eta_7\neq 0.
$$
This yields $M_8(\S^4)=\{(E\nu')\eta_7\}\cong\Z_2$. So, by Proposition \ref{rc1}(1),  $P_8(\mathbb{R}P^4)=\{\gamma_4(E\nu')\eta_7\}\cong\Z_2$. By the parallel argument, we get the group $P_9(\mathbb{R}P^4)$. We recall
$\pi_{10}(\S^4)=\{\nu^2_4,\alpha_1(4)\alpha_1(7)\}\cong\Z_{24}\oplus\Z_3$. By Lemma \ref{Y}(1) and (\ref{iRnu}),
$$
[\gamma_4\nu_4\alpha_1(7),i_\R]=\gamma_4\alpha_1(4)\alpha_1(7)=[\gamma_4\alpha_1(4)\alpha_1(7),i_\R].
$$
By (\ref{nu'nu}) and (\ref{iRnu}), we obtain $[\gamma_4\nu^2_4,i_\R]=\gamma_4\alpha_1(4)\alpha_1(7)$. This implies $\Ker[- , i_\mathbb{R}]
=\gamma_4\circ\{(\nu_4-\alpha_1(4))\alpha_1(7), 3\nu^2_4\}$, and hence, by Proposition \ref{rc1}(1), we obtain the group $P_{10}(\mathbb{R}P^4)$.

\par By Lemma \ref{PM}(3),   $M_{11}(\S^6)=\pi_{11}(\S^6)$. So, by Proposition \ref{rc1}(1), we obtain $P_{11}(\mathbb{R}P^6)
={\gamma_6}_\ast P_{11}(\mathbb{S}^6)=3\pi_{11}(\mathbb{R}P^6)$ \cite{GM}.

\par We show the case $(6,6)$. Recall $\pi_{12}(\S^6)=\{\nu^2_6\}\cong\Z_2$. By Lemma \ref{PM}(2)-(i), $M_{12}(\S^6)=\pi_{12}(\S^6)$. Hence, the assertion follows from Proposition \ref{rc1}(1).

\par Finally, recall $P_{13}(\S^6)=0$ \cite{GM}. This implies the case $(7,6)$ and completes the proof.
\end{pf}

Next, we show:
\begin{thm}\label{ex2}
Let $n\geq 2$. Then,
$P_{k+2n+1}(\mathbb{C}P^n)=
{\gamma_n}_\ast P_{k+2n+1}(\S^{2n+1})$ for $0\leq  k\leq 7$ except $(k,n)=(2,2)$ and $P_7(\C P^2)=0$, that is:

\mbox{\em (1)}
$P_{2n+1}(\mathbb{C}P^n)=\left\{ \begin{array}{ll} \pi_{2n+1}(\mathbb{C}P^n),
& \mbox{if}\ n=3,\\
2\pi_{2n+1}(\mathbb{C}P^n), & \mbox{if}\ n\neq 3;
\end{array}
\right.$

\mbox{\em (2)}
Let $k=1, 2$. Then,
$P_{k+2n+1}(\mathbb{C}P^n)=\left\{ \begin{array}{ll} 0,
& \mbox{if}\ n\ \mbox{is even},\\
\pi_{k+2n+1}(\mathbb{C}P^n), & \mbox{if}\ n \ \mbox{is odd};
\end{array}
\right.$

\mbox{\em (3)}
$P_{2n+4}(\mathbb{C}P^n)=
\left\{ \begin{array}{ll} \pi_{2n+4}(\mathbb{C}P^n),& \mbox{if}\;\;n\equiv 3\ (\bmod\,4) \ \mbox{or} \ n=2^i-2\geq 2,\\
2\pi_{2n+4}(\mathbb{C}P^n), & \mbox{if otherwise}; \end{array}\right.$

\mbox{\em (4)}
Let $k=4, 5$. Then,
$P_{k+2n+1}(\mathbb{C}P^n)=\left\{ \begin{array}{ll} 0,
& \mbox{if}\ n\geq 3,\\
\pi_{k+2n+1}(\mathbb{C}P^n), & \mbox{if}\ n=2;
\end{array}
\right.$

\mbox{\em (5)}
$P_{2n+7}(\mathbb{C}P^n)=\left\{ \begin{array}{ll} \pi_{2n+7}(\mathbb{C}P^n),& \mbox{if}\;\;n\equiv 2,3\ (\bmod\,4) \ \mbox{or} \ n=2^i-3\geq 5,\\
0, & \mbox{if otherwise}; \end{array}\right.$

\mbox{\em (6)}
$P_{2n+8}(\mathbb{C}P^n)=\left\{ \begin{array}{ll} \pi_{2n+8}(\mathbb{C}P^n),& \mbox{if}\;\;n=2,3,5 \ \mbox{or}
\ n\equiv 7\ (\bmod\,8), \\
2\pi_{2n+8}(\mathbb{C}P^n), & \mbox{if otherwise}. \end{array}\right.$
\end{thm}
\begin{pf}
Notice that, in Proposition \ref{rc1}(3), we have
$$
E^{-1}(\Ker\ {\eta_{4n+1}}_\ast)
=\left\{ \begin{array}{ll}
2\pi_{k+4n+1}(\S^{4n+1}),& \mbox{if}\;\;k=0,7,\\
0,& \mbox{if} \ k=1,2,\\
\pi_{k+4n+1}(\S^{4n+1}), & \mbox{if} \ 3\leq k\leq 6. \end{array}\right.
$$
By Lemma \ref{PM}(2)-(ii) and the relation $\eta_5\nu_6=0$, we obtain $M_{k+5,\C}(\S^5)=\pi_{k+5}(\S^5)$ for $3\leq k\leq 6$. Moreover, by the fact that $\eta_4\sigma'''=0$ \cite[(7.4)]{T} and  $[\iota_5,\eta_5]\circ Eh_0(\sigma''')=0$, we obtain $M_{12,\C}(\S^5)=\pi_{12}(\S^5)$. The assertion is obtained by these results, Proposition \ref{rc1}(1);(2), (\ref{W1})-(\ref{W6}) and the fact that $[\iota_5,\sigma''']=[\iota_7,\sigma']=0$ (see the group\ $G_{k+2n+1}(\S^{2n+1})$ for $k\leq 7$ \cite{GM}).
\end{pf}

\section{Whitehead center groups of quaternionic projective spaces}
By Lemma \ref{Y} and the fact that $\sharp[i_\H,i_\H]=12$ and $\sharp[i_\H,\gamma_n]=\frac{24}{(24,n+1)}$, we have
that $P_4(\H P^n)=[[12,\frac{24}{(24,n+1)}]]\pi_4(\H P^n)$ for $n\geq 2$, where $[[ -, - ]]$ is the least common multiple.

To state the main result, some preceding discussion is needed.
Let $p$ be an odd prime. For a generator $\alpha_{m,p}(2r+1)$ of a $p$-primary component $\pi_{k+2r+1}(\S^{2r+1};p)$, we set
$$
\alpha_{m,p}(n)=E^{n-2r-1}\alpha_{m,p}(2r+1) \ \mbox{for} \  n\geq 2r+1.
$$
There exist the elements $\alpha_{i,p}(3)\in\pi_{2i(p-1)+2}(\S^3;p)$ with $i\geq 1$ \cite[Lemma13.5]{T}. Write $\alpha_{i,3}(3)=\alpha_i(3)$ simply.
Notice that $\alpha_i(3)$ generates $\pi_{4i+2}(\S^3;3)\cong\Z_3$ for $1\leq i\leq 5$.
Recall that $\alpha_{i+1,p}(3)\in\{\alpha_{i,p}(3),p\iota_{2i(p-1)+2},\alpha_{1,p}(2i(p-1)+2)\}_1$ (Toda bracket),
$\pi_9(\S^3)=\{\alpha_1(3)\alpha_1(6)\}\\
\cong\Z_3$,\
$\pi_{10}(\S^3)=\{\alpha_2(3),\alpha_{1,5}(3)\}\cong\Z_{15}$,\  $\pi_{13}(\S^3)=\{\varepsilon',\eta_3\mu_4,\alpha_1(3)\alpha_2(6)\}\cong\Z_{12}\oplus\Z_2$,\  $\pi_{14}(\S^3)=\{\mu',\varepsilon_3\nu_{11},\nu'\varepsilon_6,\alpha_3(3),\alpha_{1,7}(3)\}\cong\Z_{84}\oplus(\Z_2)^2$\
and, by \cite[Proposition 13.6, (13.7)]{T},

\begin{equation}\label{Talpha}
\sharp(\alpha_1(n)\alpha_1(n+3))=\left\{ \begin{array}{ll}
3,&\;\mbox{if}\;\; n=3,4;\\
1,&\;\mbox{if}\;\; n\ge 5.
\end{array}
\right.
\end{equation}

Recall that there exists the element
$$
\beta_1(5)\in\{\alpha_1(5),\alpha_1(8),\alpha_1(11)\}_1\subset\pi_{15}(\mathbb{S}^5;3)\cong\Z_9.
$$
We set $\beta_1(n)=E^{n-5}\beta_1(5)$ for $n\ge 5$.
By use of \cite[Lemma 13.8, Theorem 13.9]{T} and the isomorphism $E: \pi_{15}(\S^5)\to\pi_{16}(\S^6)$,
\begin{equation}\label{12beta}
3\beta_1(5)=-\alpha_1(5)\alpha_2(8),\ \sharp\beta_1(6)=9\ \mbox{and} \ \sharp\beta_1(n)=3\ \mbox{for} \ n\geq 7.
\end{equation}

By \cite[(13.11), Theorem 13.9]{T}, $\pi_{2m+12}(\S^{2m+1};3)=\{\alpha'_3(2m+1)\}\cong\Z_9$ for $m\geq 2$ and
\begin{equation}\label{alph'3}
3\alpha'_3(5)=\alpha_3(5).
\end{equation}

We show:
\begin{prop}\label{T13.3}
 \mbox{\em(1)}
$\pi_{17}(\S^3;3)=\{\alpha_1(3)\alpha'_3(6)\}\cong\Z_3$\ and\\
$\pi_{17}(\S^3;5)=
\{\alpha_{1,5}(3)\alpha_{1,5}(10)\}\cong\Z_5$;

\mbox{\em(2)}
$\pi_{19}(\S^3;3)=\{\alpha_1(3)\alpha_1(6)\beta_1(9)\}\cong\Z_3$;

\mbox{\em(3)}
$\pi_{21}(\S^3;3)=\{\alpha_1(3)\alpha_4(6)\}\cong\Z_3$;

\mbox{\em(4)}
$\pi_{25}(\S^3;3)=\{\alpha_1(3)\alpha_5(6)\}\cong\Z_3$.
\end{prop}
\begin{pf}
We use the exact sequence \cite[Proposition 13.3]{T}. The first half of (1) follows from  (\ref{alph'3}) and \cite[Theorem 13.10.ii)]{T}. The second half of (1) follows from \cite[${\rm (13.6)}'$]{T}. By \cite[Theorem 13.10.i)]{T}, we obtain (2). (3) follows from \cite[Theorem 13.10.iii)]{T} and  (4) from \cite[Theorem 13.10.vii)]{T} and \cite[p.325]{Mi0}.
\end{pf}

We recall $\pi_k(\S^3)=\{\nu'\eta^{k-6}_6\}\cong\Z_2$ for $k=7,8$,  $\pi_{11}(\S^3)=\{\varepsilon_3\}\cong\Z_2$ and  $\pi_{12}(\S^3)=\{\eta_3\varepsilon_4,\mu_3\}\cong(\Z_2)^2$ and

\begin{equation}\label{nu'nu}
\nu'\nu_6=0.
\end{equation}

We also recall:
$\pi_{15}(\S^3)=\{\nu'\mu_6,\nu'\eta_6\varepsilon_7\}\cong(\Z_2)^2$,\  $\pi_{16}(\S^3)=\\
\{\nu'\eta_6\mu_7,\alpha_1(3)\beta_1(6)\}\cong\Z_6$,\ $\pi_{17}(\S^3)=
\{\varepsilon_3\nu^2_{11},\alpha_1(3)\alpha'_3(6),\alpha_{1,5}(3)\alpha_{1,5}(10)\}\cong\Z_{30}$\ and\
$\pi_{18}(\S^3)=\{\bar{\varepsilon}_3,\alpha_4(3),\alpha_{2,5}(3)\}\cong\Z_{30}$.

Now, we show:
\begin{lem}\label{Q}
\mbox{\em (1)}
$Q_3(\S^3)=12\pi_3(\S^3)$\ and\
$Q_k(\S^3)=0$ for $k=4,5,11,12$;

\mbox{\em (2)}
$Q_6(\S^3)=3\pi_6(\S^3)\cong\Z_4, Q_k(\S^3)=\pi_k(\S^3)\cong\Z_2$ for $k=7,8$ and  $Q_9(\S^3)=\pi_9(\S^3)\cong\Z_3$;

\mbox{\em (3)}
$Q_{10}(\S^3)=3\pi_{10}(\S^3)\cong\Z_5$\ and\  $Q_{13}(\S^3)=\{\varepsilon',\alpha_1(3)\alpha_2(6)\}\cong\Z_{12}$;

\mbox{\em (4)}
$Q_{14}(\S^3)=
\pi_{14}(\S^3)\cong\Z_{84}\oplus(\Z_2)^2$;

\mbox{\em (5)}
$Q_{15}(\S^3)=\pi_{15}(\S^3)\cong(\Z_2)^2$,\  $Q_{16}(\S^3)=\pi^3_{16}=\{\nu'\eta_6\mu_7\}\cong\Z_2$,\
$Q_{17}(\S^3)=\pi_{17}(\S^3)\cong\Z_{30}$\ and\
$Q_{18}(\S^3)=6\pi_{18}(\S^3)\cong\Z_5$.
\end{lem}
\begin{pf}
Since $\sharp\omega=12$, we have the first of (1).
We know $\nu'E^3\alpha\neq 0$ for  $\alpha=\eta_3,\eta^2_3,\varepsilon_3,\eta_3\varepsilon_4,\mu_3$ and $\eta_3\mu_4$. This leads to the rest of (1) and the fact that $Q_{13}(\S^3)\not\ni\eta_3\mu_4$. By (\ref{Talpha}) and (\ref{nu'nu}),  $\omega E^3\omega=E(\alpha_1(3)\alpha_1(6))$. Hence,  $Q_6(\S^3)=3\pi_6(\S^3)\cong\Z_4$.
Obviously, $Q_{k+6}(\S^3)=\pi_{k+6}(\S^3)$ for $k=1,2,3$ and (2) follows.

By (\ref{12beta}),
$\alpha_2(3)\not\in Q_{10}(\S^3)$ and
$Q_{10}(\S^3)=\{\alpha_{1,5}(3)\}\cong\Z_5$.

We recall
$\varepsilon'\in\{\nu',2\nu_6,\nu_9\}=\{\nu',\nu_6,2\nu_9\}$. By \cite[Lemma 6.6,(7.12)]{T},
$$
\nu'E^3\varepsilon'\in\nu'\circ\{2\nu_6,\nu_9,2\nu_{12}\}=\{\nu',2\nu_6,\nu_9\}\circ 2\nu_{13}=\varepsilon'\circ 2\nu_{13}=2\eta^2_3\varepsilon_5\nu_{13}=0.
$$
By
(\ref{12beta}), $Q_{13}(\S^3)\ni\alpha_1(3)\alpha_2(6)$.
This leads to (3).

By \cite[Proposition (2.2)(5)]{Og}, we know
$\eta_5\zeta_6
=0.$
So, by the relation $2\zeta_5=\pm E^2\mu'$ \cite[(7.14)]{T}, we obtain
$\nu'E^3\mu'=\eta^3_3\zeta_6=0$. Hence, $\mu'\in Q_{14}(\S^3)$. Obviously, $\nu'\varepsilon_6,\alpha_{1,7}(3)\in Q_{14}(\S^3)$.
By the relation \cite[(7.12)]{T}:
\begin{equation}\label{t712}
\varepsilon_3\nu_{11}=\nu'\bar{\nu}_6,
\end{equation}
we have $\nu'\varepsilon_6\nu_{14}=0$ and $
\varepsilon_3\nu_{11}\in Q_{14}(\S^3)$.
The fact that $\alpha_3(3)\in Q_{14}(\S^3)$ follows from the relation
$\alpha_1(3)\alpha_3(6)=3\alpha_1(3)\alpha'_3(6)=0$ (\ref{alph'3}), and (4) follows.

The first of (5) follows from
the group structure
$\pi_{15}(\S^3)=\{\nu'\mu_6,\nu'\eta_6\varepsilon_7\}\cong(\Z_2)^2$, the fact that $\nu'\in Q_6(\S^3)$ (2)
 and Lemma \ref{LQ}(1). By  the relation $\nu'\eta_6\mu_7=\nu'\mu_6\eta_{15}$ and Lemma \ref{LQ}(1),
 $\nu'\eta_6\mu_7\in Q_{16}(\S^3)$. So, by Proposition \ref{T13.3}(2), we obtain
$Q_{16}(\S^3)=\pi^3_{16}=\{\nu'\eta_6\mu_7\}\cong\Z_2$ and leads to the second of (5).

By the relation (\ref{t712}), we have $\nu'(\varepsilon_6\nu^2_{14})=0$. In the exact sequence ($p=3,i=19$ in \cite[Proposition 13.3]{T}):
$$
\pi_{19}(\S^5;3)\rarrow{G}\pi_{20}(\S^3;3)\rarrow{H}\pi_{20}(\S^7;3)\rarrow{\Delta}\pi_{18}(\S^5;3),
$$
we know that all groups are isomorphic to $Z_3$ \cite[Theorem 13.10]{T} and that $\Delta(E^2(\alpha_1(5)\beta_1(8))=3(\alpha_1(5)\beta_1(8))=0$ and $G(\beta)=\alpha_1(3)E\beta$ for $\beta\in\pi_{19}(\S^5;3)$. This implies the relation $\alpha_1(3)E\beta=0$ and this yields the third of (5).

We know  $\nu'\bar{\varepsilon}_6\neq 0$ \cite[Theorem 12.8]{T} and $\alpha_1(3)\alpha_4(6)\neq 0$ (Proposition \ref{T13.3}(3)).  This leads to the last of (5) and completes the proof.
\end{pf}

Set $\sigma',\sigma'',\sigma'''$ for the generators of $\pi^n_{n+7}\cong\Z_{2^{n-4}}$ with $5\leq n\leq 7$.
By abuse of notation, these are also used to represent the generators of $\pi_{n+7}(\mathbb{S}^n)\cong\Z_{15\cdot 2^{n-4}}$ with $5\leq n\leq 7$, respectively. We show:
\begin{prop}\label{nu6P}
Let $4\leq  k\leq 22$. Then,
$\nu_6E^6\pi^3_k\subseteq P\pi^{13}_{k+8}$, $\nu_7E^7\pi_k(\S^3)=0$ for $4\leq  k\leq 21$ and $\nu_kE^k\pi_{22}(\S^3)\cong\Z_3$ for $5\leq k\leq 9$.
\end{prop}
\begin{pf}
First of all, recall that
$\{\nu_6,\eta_9,2\iota_{10}\}=[\iota_6,\iota_6]+\{2[\iota_6,\iota_6]\}$.
We know $\nu_5\circ E^5\nu'=2\nu^2_5=0$. This implies
$\nu_6E^6\pi^3_k=0$ for $4\leq  k\leq 10$.

We know $[\iota_6,\eta_6]=0$. So, we obtain $[\iota_6,\eta_6\varepsilon_7]=0$ and $[\iota_6,\nu^3_6]=[\iota_6,\eta_6\bar{\nu}_7]=0$. Hence, by \cite[(10.12)]{T}, $[\iota_6,\mu_6]\equiv 0\ \bmod\ [\iota_6,\eta_6\varepsilon_7], [\iota_6,\nu^3_6]$.
That is
\begin{equation}\label{Whi6mu}
[\iota_6,\mu_6]=0.
\end{equation}

Let $\alpha\in\pi^5_k$ and $\beta\in\pi^3_k$ be elements such that $2\iota_5\circ\alpha=0$ and  $\{\eta_4,2\iota_5,\alpha\}\ni E\beta$. Then,
$$
\nu_6E^6\beta\in\nu_6\circ\{\eta_9,2\iota_{10},E^5\alpha\}=[\iota_6,\iota_6]\circ E^6\alpha.
$$
By \cite[(7.25)]{T},

\begin{equation}\label{wsg''}
\nu_6\mu_9=8P(\sigma_{13})=2[\iota_6,\sigma'']
=[\iota_6,E\sigma'''].
\end{equation}
This and (\ref{nu6ep}) imply
$\nu_6E^6\pi^3_{11}=P\pi^{13}_{19}$ and $\nu_6E^6\pi^3_{12}=8P\pi^{13}_{20}$. From the fact that $\mu'\in\{\eta_3,2\iota_4,\mu_4\}$ and (\ref{Whi6mu}), we have
\begin{equation}\label{Whmu6}
\nu_6E^6\mu'\in\nu_6\circ\{\eta_9,2\iota_{10},\mu_{10}\}=[\iota_6,\mu_6]=0.
\end{equation}

By the fact that $\pi^5_{18}\cong(\Z_2)^2$, we obtain $\nu_5E^5\varepsilon'\in \{\nu_5,2\nu_8,\nu_{11}\}\circ 2\nu_{15}\subseteq
2\pi^5_{18}=0$. We obtain $\nu_6E^6(\nu'\varepsilon_6)=0$.
By the relation \cite[(7.18)]{T}:
\begin{equation}\label{Whnu36}
2\bar{\nu}_6\circ\nu_{14}=\nu_6\varepsilon_9=[\iota_6,\nu^2_6],
\end{equation}
we have $\nu_6\varepsilon_9\nu_{17}=2\bar{\nu}_6\circ\nu^2_{14}=0$.

By (\ref{nu'nu}), (\ref{t712}) and the relation $\nu'\zeta_6=0$ \cite[Proposition (2.4)]{Og}, we have
$$
\pi^3_k=\nu'\circ\pi^6_k,\ 15\le k\le 17.
$$
Therefore, we obtain $\nu_6E^6\pi^3_k=0$ for $15\le k\le 17$.

Since $\pi^3_{18}=\{\bar{\varepsilon}_3\}\cong\Z_2$ and $\bar{\varepsilon}_6=\eta_6\kappa_7$ \cite[(10.23)]{T},
we have $\nu_6E^6\pi^3_{18}=0$.

Since $\pi^3_{19}=\{\eta_3\bar{\varepsilon}_4,\mu_3\sigma_{12}\}\cong(\Z_2)^2$, $\nu_6\eta_9\bar{\varepsilon}_{10}=0$ and $\nu_6\mu_9\sigma_{18}=8P(\sigma^2_{13})$ (\ref{wsg''}), we have $\nu_6E^6\pi^3_{18}=8\{P(\sigma^2_{13})\}$.

Next we recall
$$
\pi^3_{20}=\{\bar{\varepsilon}',\bar{\mu}_3,\eta_3\mu_4\sigma_{12}\}\cong\Z_4\oplus(\Z_2)^2
$$
and $E\bar{\varepsilon}'=(E\nu')\kappa_7$ \cite[Lemma 12.3]{T}. By \cite[(16.6)]{MT}, $\nu_6\bar{\mu}_9=16P(\rho_{13})$. This implies $\nu_6E^6\pi^3_{20}=16\{P(\rho_{13})\}$.

From the fact that $\pi^3_{21}=\{\mu'\sigma_{14},\nu'\bar{\varepsilon}_6,\eta_3\bar{\mu}_4\}$ and (\ref{Whmu6}), we obtain
$\nu_6E^6\pi^3_{21}=0$.

We recall $\pi^6_{21}\cong\Z_4\oplus\Z_2$.
Since $\pi^3_{21}=\{\bar{\mu}',\nu'\mu_6\sigma_{15}\}\cong\Z_4\oplus\Z_2$, by (\ref{Whmu6}), the definition of $\mu'$ and \cite[(3.2)]{MMO}, we see that
$$
\nu_6E^6\bar{\mu}'\in\nu_6\{E^6\mu',4\iota_{20},\sigma_{20}\}\subset\{0,4\iota_{20},4\sigma_{20}\}=
\pi^6_{21}\circ 4\sigma_{21}=0.
$$
This implies $\nu_6E^6\pi^3_{22}=0$.

Next, we examine the case for odd primes. This case appears for $k=6,9,10,13,14$ and $16\leq k\leq 22$. By (\ref{Talpha}), the assertion
$$
(\star) \hspace{1cm} \alpha_1(7)E^7\pi_k(\S^3;3)=0
$$
holds for $k=6,9,13,16$.

By (\ref{12beta}) and (\ref{alph'3}),
$(\star)$ holds for $k=10,14$.

By (\ref{Talpha}) and Proposition \ref{T13.3}(1), $(\star)$ holds for $k=17$.

For $k=18$, we have $\alpha_1(7)E^7\pi_{18}(\S^3;3)=\{\alpha_1(7)\alpha_4(10)\}\subseteq E^2\pi_{23}(\S^5;3)=0$  \cite[p.185]{T}.

By \cite[Theorem 13.10.iv);v)]{T}, $E^2\pi_{l+19}(\S^3;3)=0$ for $l=0,1$, and hence, $(\star)$ holds for $k=19,20$.

By \cite[p.185]{T}, $E^2\pi_{21}(\S^3;3)=0$, and hence, $(\star)$ for $k=21$ holds.

Finally, for $k=22$, by Proposition \ref{T13.3}(4) and  \cite[p.56]{T2}, $\alpha_1(9)\alpha_5(12)\neq 0$. This leads to the last result and completes the proof.
\end{pf}

We propose:\\
\bigskip

\noindent
{\bf Problem.}\ If $\nu'E^3\beta=0$ for $\beta\in\pi^3_{k-1}$, then $\nu_4\wedge\beta=0$?
\bigskip

Lemma \ref{LQ}(2) and Proposition \ref{nu6P} yield:
\begin{cor}\label{L''}
$L''_{k-1,n}(\S^3)=\pi_{k-1}(\S^3)$\ for\ $n\geq 2$ and $5\leq k\leq 22$.
\end{cor}

Writing $R_k^n=\pi_k(SO(n);2)$, we get from $({\mathcal{SF}}^n_k)$ the exact sequence
$$
({\mathcal{R}}^n_k) \ \cdots\rarrow{}
\pi_{k+1}^n\rarrow{\Delta}R_k^n\rarrow{i_\ast} R_k^{n+1}\rarrow{p_\ast}\cdots.
$$

Next, we show:
\begin{thm}\label{PHP1}
$P_{4n+3+k}(\H P^n)=P'_{4n+3+k}(\H P^n)\oplus P''_{4n+3+k}(\H P^n)$ for $0\le k\le 10$, where

{\em (1)}
$P'_{4n+3}(\mathbb{H}P^n)=[[\frac{24}{(24,n+1)},2]]{\gamma_n}_\ast\pi_{4n+3}(\S^{4n+3})$ for $n\ge 2$;

{\em (2)}
$P'_{4n+3+k}(\mathbb{H}P^n)={\gamma_n}_\ast\pi_{4n+3+k}(\S^{4n+3})$ for $k=1,2,4,5,8,9,10$;

{\em (3)}
$P'_{4n+6}(\mathbb{H}P^n)=\frac{3+(-1)^n}{2}{\gamma_n}_\ast\pi_{4n+6}(\S^{4n+3})$;

{\em (4)}
$P'_{4n+9}(\mathbb{H}P^n)=
\frac{1-(-1)^n}{2}{\gamma_n}_\ast\pi_{4n+9}(\mathbb{S}^{4n+3})$;

{\em (5)} $P'_{4n+10}(\mathbb{H}P^n)=\left\{ \begin{array}{ll} 2{\gamma_n}_\ast\pi_{4n+10}(\S^{4n+3}),
& \mbox{if}\ n\equiv 0,1,2\ (\bmod\ 4);\\
{\gamma_n}_\ast\pi_{4n+10}(\S^{4n+3}), & \mbox{if}\ n\equiv 3\ (\bmod\ 4) \ \mbox{or} \ n=2.
\end{array}
\right.$
\end{thm}
\begin{pf}
In view of Lemma \ref{LQ}(2),(3), Proposition \ref{newP} and Proposition \ref{rc12}(1), we get
$$
P_k(\mathbb{H}P^{2n+1})= P'_k(\mathbb{H}P^{2n+1})\oplus P''_k(\mathbb{H}P^{2n+1}).
$$
We show that
$$
P_{8n+3+k}(\mathbb{H}P^{2n})= P'_{8n+3+k}(\mathbb{H}P^{2n})\oplus P''_{8n+3+k}(\mathbb{H}P^{2n})
$$
for $0\le k\le 10$.
In virtue of (\ref{W1})-(\ref{W6}), we can determine
$P_{4n+3+k}(\S^{4n+3})$ for $0\le k\le 7$. Hence, by Proposition \ref{newP} and Lemma \ref{PM}(8)-(11), we obtain $P'_{4n+3+k}(\H P^n)$ for $0\le k\le 7$. In the case that $8\le k\le 10$, we have $M_{4n+3+k}(\S^{4n+3})=G_{4n+3+k}(\S^{4n+3})=\pi_{4n+3+k}(\S^{4n+3})$. This determines $P'_{4n+3+k}(\H P^n)$ for $0\le k\le 10$.

Now our task is to check the condition $(\heartsuit)$ in Proposition \ref{rc12}.2(iii). Let $n$ be even and we work in the $2$-primary component. The assertion is a direct consequence of the fact that $[\iota_{4n+3},\pi^{4n+3}_{4n+3+k}]=0$ for $k=1,2,8,9,10$ \cite[Proposition 1.3]{GM}.

For $k=0,3,6$, the condition $(\heartsuit)$ holds because $2[\iota_{4n+3},\iota_{4n+3}]=2[\iota_{4n+3},\nu_{4n+3}]=2[\iota_{4n+3},\nu^2_{4n+3}]=0$.

We know
$[\iota_{4n+3},\sigma_{4n+3}]=0$ if and only if $n=2$ or $n\equiv 3\ (\bmod\ 4)$ (\ref{W6}). So, for the case $n\ne 2$ and $n\equiv 0,2\ (\bmod\ 4)$, suppose that
there exists an element $\beta\in\pi^3_{4n+9}$ such that $J(\Delta\sigma_{4n+3})=[\iota_{4n+3},\sigma_{4n+3}]=\nu_{4n+3}E^{4n+3}\beta=J([\iota_3]_{4n+3}\beta)$.
We set $R^n_k=\pi_k(SO(n);2)$. By \cite{BM} and \cite{HM},
$$
R_{4n+9}^{4n+k}\cong\left\{\begin{array}{ll}
\Z_2\oplus\Z_8\oplus\Z_2,&\;\mbox{if}\;\; k=1;\\
\Z_2\oplus\Z_4,&\;\mbox{if}\;\; k=2;\\
\Z_2\oplus\Z_2,&\;\mbox{if}\;\; k=3;\\
\Z_2,&\;\mbox{if}\;\;k=4.
\end{array}
\right.
$$
Here,
the direct summand $\Z_2$ corresponds to the Bott result: $R^{\infty}_{4n+9}\cong\Z_2$.
Since $[\iota_{4n+3},\sigma_{4n+3}]\ne 0$ and this element does not correspond to the Bott element,
$J: R^{4n+3}_{4n+9}\to\pi^{4n+3}_{8n+12}$ is a monomorphism and hence,  $\Delta\sigma_{4n+3}=[\iota_3]_{4n+3}\beta$.
In the exact sequence $({\mathcal{R}}^{4n+k-1}_{4n+9})\ (k=2,3)$:
$$
R^{4n+k-1}_{4n+9}\rarrow{i_*}R^{4n+k}_{4n+9}
\rarrow{p_*}\pi^{4n+k-1}_{4n+9},
$$
we know $\sharp[\iota_{4n+2},\sigma_{4n+2}]=16$ (\ref{W6}) and $
[\iota_{4n+1},\eta_{4n+1}\sigma_{4n+2}]\ne 0,\ [\iota_{4n+1},\bar{\nu}_{4n+1}]\ne 0$ and $[\iota_{4n+1},\varepsilon_{4n+1}]\ne 0$ \cite[Lemma 4.3]{GM}. So, $i_*: R^{4n+k-1}_{4n+9}\to R^{4n+k}_{4n+9}$ are epimorphisms for $k=2,3$, respectively. This shows that $[\iota_3]_{4n+1}\beta$ generates the direct summand $\Z_8$
in $R^{4n+1}_{4n+9}$ and contradicts the fact that $4\beta=0$ \cite[Corollary (1.22)]{J2}. Hence, we have the assertion for $k=7$. This completes the proof.
\end{pf}

For $n$ even,
since $[\iota_{8n+3},\sigma^2_{8n+3}]\ne 0 $ \cite[Table 1]{Mu}, The assertion of Theorem \ref{PHP1} for $k=7$, $n\equiv \ 0\  (\bmod\ 4)$ is obtained from Corollary \ref{alphal}(2).

Since $[\iota_n,\zeta_n]=0$ except for $n\not\equiv 115\ (\bmod\ 128)$ \cite[p. 426]{GM}, we obtain

\begin{rem} {\em
The assertion of Theorem \ref{PHP1} for $k=11$ holds and $P'_{4n+14}(\H P^n)={\gamma_n}_*\pi_{4n+14}(\S^{4n+3})$ if
$n\not\equiv 115\ (\bmod\ 128)$.}
\end{rem}

By Proposition \ref{rc12}, Lemma \ref{Q} and
Theorem \ref{PHP1}, we obtain:

\begin{prop}\label{exh1}
\mbox{\em (1)}
Let $n\geq 2$. Then:

$
P_k(\mathbb{H}P^n)=\left\{ \begin{array}{ll}
0, & \mbox{if} \hspace{3mm} k=5,6,\\
{i_{\mathbb{H}}}_\ast E\{\nu'\}\cong\Z_4,
& \mbox{if} \hspace{3mm} k=7,\\
{i_{\mathbb{H}}}_\ast E\{\nu'\eta_6\}\cong\Z_2, & \mbox{if} \hspace{3mm} k=8,\\
{i_{\mathbb{H}}}_\ast E\{\nu'\eta^2_6\}\cong\Z_2,&\mbox{if} \hspace{3mm} k=9,\\
{i_{\mathbb{H}}}_\ast E\{\alpha_1(3)\alpha_1(6)\}\cong\Z_3,&\mbox{if} \hspace{3mm} k=10;
\end{array}
\right.
$
\vspace{1mm}

\mbox{\em (2)}
$P_{11}(\mathbb{H}P^2)=\{8\gamma_2,i_\H E\alpha_{1,5}(3)\}\cong 8\Z\oplus\Z_5$;

\mbox{\em (3)}
$P_k(\mathbb{H}P^2)={\gamma_2}_\ast \pi_k(\S^{11})\cong\Z_2$, if \ $k=12,13$;

\mbox{\em (4)}
Let $n\geq 3$. Then,

$P_k(\mathbb{H}P^n)=\left\{ \begin{array}{ll}{i_{\mathbb{H}}}_\ast E\{\alpha_{1,5}(3)\}\cong\Z_5,
& \mbox{if} \hspace{3mm} k=11,\\
0,&\mbox{if} \hspace{3mm} k=12,13,\\
{i_{\mathbb{H}}}_\ast E\{\varepsilon',\alpha_1(3)\alpha_2(6)\}\cong\Z_{12},&\mbox{if} \hspace{3mm} k=14;
\end{array}
\right.
$

\vspace{1mm}

\mbox{\em (5)}
$P_{15}(\mathbb{H}P^3)=\{6\gamma_3\}\oplus
{i_{\mathbb{H}}}_\ast E\pi_{14}(\S^3)
\cong 6\Z\oplus\Z_{84}\oplus(\Z_2)^2$;

\mbox{\em (6)}
$P_{16}(\mathbb{H}P^3)=\pi_{16}(\mathbb{H}P^3)\cong(\Z_2)^3$,\ $P_{17}(\mathbb{H}P^3)={\gamma_3}_\ast\pi_{17}(\mathbb{S}^{15})\oplus{i_\H}_\ast E\pi^3_{16}
\cong(\Z_2)^2$\ and \ $P_{18}(\mathbb{H}P^3)=\pi_{18}(\mathbb{H}P^3)\cong\Z_{24}\oplus\Z_{30}$;

\end{prop}

\section{Gottlieb groups of real projective spaces}
Hereafter, we set
$$
G'_k(\mathbb{F}P^n)=G_k(\mathbb{F}P^n)\cap({\gamma_n}_\ast\pi_k(\S^{d(n+1)-1})).
$$
Notice that $G'_k(\F P^n)=G_k(\F P^n)$ if $\F=\R,\C$ and $k\ge d+1$.

First of all, we show

\begin{prop}\label{Gotc}
$G'_k(\F P^n)\subseteq{\gamma_n}_*G_k(S^{d(n+1)-1})$.
\end{prop}
\begin{pf}
The real case is just \cite[Theorem 6-2]{G1}.
For any element $\gamma_n\alpha\in G'_k(\F P^n)$, we have $0=[\gamma_n\alpha,\gamma_n]=\gamma_n\circ[\alpha,\iota_{d(n+1)-1}]$. Since\\
${\gamma_n}_*: \pi_{k+d(n+1)-2}(S^{d(n+1)-1})\to\pi_{k+d(n+1)-2}(\F P^n)$ is a monomorphism, $\alpha\in P_k(S^{d(n+1)-1})=G_k(S^{d(n+1)-1})$. 
\end{pf}

In the sequel, we need the following. Let $K$ be a closed subgroup of a Lie group $H$ and write $H/K$ for the left coset.
We recall \cite[Example 2.2]{Si}, \cite[Theorem II.5]{L}, which is directly obtained by Lemma \ref{SV} and
the fact that the usual pairing $H\times H/K\to H/K$ is an associated map with the projection $p: H\to H/K$:

\begin{lem} \label{SL}
Let $K$ be a closed subgroup of a Lie group $H$ and $H/K$ the left coset. Then, the projection $p: H\to H/K$ is cyclic and $p_*(\pi_n(H))\subseteq G_n(H/K)$ for $n\geq 1$.
\end{lem}

Write $i'_{n,\mathbb{F}}: O_\mathbb{F}(n-1)\times O_\mathbb{F}(1)\hookrightarrow O_\mathbb{F}(n)$ for the inclusion,
$p'_{n,\mathbb{F}}: O_\mathbb{F}(n)\to\mathbb{F}P^{n-1}$ for the quotient map.
Now, we consider the exact sequence induced from the fibration $O_{\mathbb{F}}(n+1)\stackrel{O_{\mathbb{F}}(n)\times O_{\mathbb{F}}(1)}{\longrightarrow}\mathbb{F}P^n$:
$$
(\mathcal{FP}^n_k) \ \cdots\rarrow{}
\pi_k(\mathbb{F}P^n)\rarrow{\Delta'_\mathbb{F}}\pi_{k-1}(O_{\mathbb{F}}(n)\times O_{\mathbb{F}}(1))\rarrow{i'_\ast}\pi_{k-1}(O_{\mathbb{F}}(n+1))\rarrow{p'_\ast}\cdots,
$$
where $i'=i'_{n+1,\mathbb{F}}$ and $p'=p'_{n+1,\mathbb{F}}$.
Then, by Lemma \ref{SL}, we have
\begin{equation}\label{sl}
\Ker\ \Delta'_\F=p'_*\pi_k(O_\mathbb{F}(n+1))\subseteq G_k(\K P^n).
\end{equation}

\bigskip

\par Next, we consider the natural map from $(\mathcal{SF}^n_k)$ to $(\mathcal{FP}^n_k)$:

\begin{equation}\label{comg}
\xymatrix{
\pi_k(SO_\mathbb{F}(n+1))\ar@{=}[d]\ar[r]^(.6){p_\ast}&
\pi_k(\mathbb{S}^{d(n+1)-1})\ar[d]^{{\gamma_n}_\ast}\ar[r]^(.4){\Delta_\mathbb{F}}&\pi_{k-1}(SO_\mathbb{F}(n))\ar^{\cap}[d]\\
\pi_k(O_\mathbb{F}(n+1))\ar[r]^{p'_\ast}&\pi_k(\mathbb{F}P^n)\ar[r]^(.4){\Delta'_\mathbb{F}}&\pi_{k-1}(O_\mathbb{F}(n)\times O_\mathbb{F}(1)).}
\end{equation}

We show:
\begin{lem}\label{fund1}
\mbox{\em (1)} $\Ker\{\Delta_\mathbb{F}: \pi_k(\mathbb{S}^{d(n+1)-1})\to\pi_{k-1}(SO_\mathbb{F}(n))\}
\subseteq{{\gamma_n}}^{-1}_\ast G_k(\mathbb{F}P^n)$.

\mbox{\em (2)} Let $k\geq d+1$. If $E^{d-1}\circ J_\mathbb{F}\mid_{\Delta_\mathbb{F}(\pi_k
(\mathbb{S}^{d(n+1)-1}))} : \Delta_\mathbb{F}(\pi_k(\mathbb{S}^{d(n+1)-1}))\to \pi_{k+d(n+1)-2}(\mathbb{S}^{d(n+1)-1})$ is a monomorphism, then $\gamma_\ast G_k(\mathbb{S}^{d(n+1)-1})\subseteq G_k(\mathbb{F}P^n)$.
In particular, under the assumption,\\
$G_k(\mathbb{F}P^n)={\gamma_n}_\ast G_k(S^{d(n+1)-1})$ for $\F=\R, \C$ and\  $G'_k(\mathbb{H}P^n)={\gamma_n}_\ast G_k(S^{4n+3})$.
\end{lem}
\begin{pf} By the commutativity of the right square of the diagram (\ref{comg}),
$\Ker\ \Delta_\mathbb{F}={{\gamma_n}}^{-1}_\ast(\Ker\ \Delta'_\mathbb{F})$. Hence, (\ref{sl}) implies (1).

By the assumption,
(\ref{JDel}) and (\ref{JDel2}), $G_k(\mathbb{S}^{d(n+1)-1})=\Ker(J\circ\Delta)=\Ker\ \Delta_\mathbb{F}$. So, by (1), we obtain
$\gamma_\ast G_k(\mathbb{S}^{d(n+1)-1})\subseteq G_k(\mathbb{F}P^n)$. By Proposition \ref{Gotc},
$G_k(\mathbb{F}P^n)\subseteq 
\gamma_\ast G_k(\mathbb{S}^{d(n+1)-1})$ for $\F=\R, \C$ and $G'_k(\mathbb{F}P^n)\subseteq \gamma_\ast G_k(\mathbb{S}^{d(n+1)-1})$ for $\F=\H$.
This completes the proof.
\end{pf}

By Proposition \ref{Gotc} and Lemma \ref{fund1}(1), we obtain:
\begin{cor}\label{GRP}
If $\Delta_\R\pi_k(\S^n)=0$, then $G_k(\R P^n)={\gamma_n}_*G_k(\S^n)$.
\end{cor}

By \cite{G} and \cite{PW1}, we know:
\begin{thm}{\bf (Gottlieb)} \label{GO}
$G_1(\mathbb{R}P^n)=\left\{ \begin{array}{ll} 0,
& \mbox{if}\;\;n\;\;\mbox{is even};\\
\pi_1(\mathbb{R}P^n), & \mbox{if}\;\;n\;\;\mbox{is odd.}
\end{array}
\right.$
\end{thm}

\begin{thm}{\bf (Pak-Woo)} \label{PW}

$G_n(\mathbb{R}P^n)={\gamma_n}_\ast G_n(\mathbb{S}^n)=\left\{ \begin{array}{ll} 0,
& \mbox{if}\;\;n\;\;\mbox{is even};\\
\pi_n(\mathbb{R}P^n), & \mbox{if}\;\;n=1,3,7;\\
2\pi_n(\mathbb{R}P^n), & \mbox{if}\;\;n\;\;\mbox{is odd}\;\;\mbox{and}\;\;n\not=1,3,7.
\end{array}
\right.$
\end{thm}

By Lemma \ref{fund1}(1) and the fact that $\pi_{k-1}(SO(2))=0$ for $k\geq 3$, we derive:
\begin{prop}
$G_k(\mathbb{R}P^2)=\pi_k(\mathbb{R}P^2)$ for $k\ge 3$.
\end{prop}

Theorems \ref{main0}, \ref{GO} and \ref{PW} imply $G_k(\mathbb{R}P^n)=P_k(\mathbb{R}P^n)$ for $k=1,n$.

\bigskip

We examine the the equality:
$$
(\ast)_\alpha \hspace{1cm}
\sharp\Delta\alpha=\sharp[\iota_n,\alpha] \ \mbox{for} \ \alpha\in\pi_{k+n}(\S^n).
$$

Let $p$ be an odd prime.
By use of Serre's isomorphism \cite[(13.1)]{T}, the equality\  $(\ast)_{E\alpha}$\ for\ $\alpha\in\pi_{k+2n-1}(\S^{2n-1};p)$\ holds if\ $\sharp\alpha=\sharp E^{2n}\alpha$. For example,
$\alpha$ is taken as $\alpha_1(2m+1)\in\pi_{2m+4}(\S^{2m+1};3),
\alpha_2(2m+1)\in\pi_{2m+8}(\S^{2m+1};3), \alpha_{1,5}(2m+1)\in\pi_{2m+8}(\S^{2m+1};5)$ for $m\geq 1$ and $\beta_1(2m+1)\in\pi_{2m+11}(\S^{2m+1};3)$ for $m\geq 3$. We show:
\begin{lem}\label{nlem}
Let $\alpha\in\pi_{n+k}(\mathbb{S}^n)$ for $k\leq 7$.
Then:

\mbox{\em (1)}\
The inequality $\sharp\Delta\alpha\gneqq\sharp[\iota_n,\alpha]$
holds if and only if: $\alpha=E\nu'$, $(E\nu')\eta_7$, $\nu_4\eta_7,(E\nu')\eta^2_7$,
$\nu_4\eta^2_7$, $[\iota_6,\iota_6]$, $E\sigma'$, $\sigma_8$, $\sigma_{11}$,
$\nu_n$ $\mbox{with}\ n=2^i-3\; \mbox{for}\;i\geq 4$ or $\nu^2_n\, \mbox{with}\ n=2^i-5\;\mbox{for}\;i\geq 4$;

\mbox{\em (2)}\
The restriction map
$J\mid_{\Delta\pi_k(\mathbb{S}^n)} : \Delta\pi_k(\mathbb{S}^n)\to\pi_{k+n-1}(\mathbb{S}^n)$ is not a monomorphism if and only if $(k,n)=(3,4),(4,4),(5,4),(6,4),(5,6),(7,8),\\
(7,11),(3,2^i-3)\;
\mbox{with}\; i\ge 4\;\mbox{and}\;(6,2^i-5)\;\mbox{with $i\ge 5$}$.

\end{lem}
\begin{pf}
By use of the exact sequence $({\mathcal R}^{n}_{n+k})$ for $-1\leq k\leq 1$ and \cite[pp. 161-2, Tables]{K}, we obtain
$$
\sharp\Delta\iota_n=\left\{\begin{array}{ll}
1,&\;\mbox{if}\ n=1,3,7;\\
2,&\;\mbox{if $n$ is odd and}\ n\neq 1,3,7;\\
\infty,&\;\mbox{if $n$ is even};\\
\end{array}
\right.
$$
$$
\sharp\Delta\eta_n=\left\{\begin{array}{ll}
1,&\;\mbox{if}\ n=2,6\ \mbox{or}\ n\equiv 3\ (\bmod\ 4);\\
2,&\;\mbox{if otherwise};
\end{array}
\right.
$$
$$
\sharp\Delta\eta^2_n=\left\{\begin{array}{ll}
1,&\;\mbox{if}\ n\equiv 2,3\ (\bmod\ 4);\\
2,&\;\mbox{if otherwise}.
\end{array}
\right.
$$
By (\ref{JDel2}), (\ref{W1}), (\ref{W2}) and (\ref{W3}), the equality $(\ast)_\alpha$ holds for $\alpha=\iota_n,\eta_n,\eta^2_n$.

We show $\sharp\Delta(E\nu')=4$.
In the exact sequence $(\mathcal{R}^4_3)$:
$$
\pi_4(\mathbb{S}^4)
\stackrel{\Delta}{\to}\pi_3(SO(4))\stackrel{i_\ast}{\to}\pi_3(SO(5))\stackrel{}{\to}0,
$$
by \cite[Theorem 23.6]{St},
\begin{equation}\label{SO4}
\pm\Delta\iota_4=2[\iota_3]-[\eta_2]_4.
\end{equation}
This implies
$\pm\Delta(E\nu')=\pm\Delta\iota_4\circ\nu'=2[\iota_3]\nu'-[\eta_2]_4\nu'$ and hence $\sharp\Delta(E\nu')=4$.
On the other hand
$\sharp[\iota_4,E\nu']=2$ (\ref{W40}).
We recall
\begin{equation}\label{7R5}
\pi_7(SO(5))\cong\Z,
\end{equation}
$\pi_k(SO(5))=0$ for $k=8,9$ and
$J(\Delta\nu_4)=\pm[\iota_4,\nu_4]=\pm 2\nu^2_4$ by Lemma \ref{Y}(3)
and (\ref{JDel2}). So, in the sequence $(\mathcal{R}^4_6)$, we
get that
\begin{equation}\label{Dnu4}
\Delta\nu_4\equiv\pm[\iota_3]\omega\ (\bmod\ 3[\eta_2]_4\omega)
\end{equation}
and that the equality $(\ast)_{\nu_4}$ holds.
\par By (\ref{7R5}), we have $\Delta\nu_5=0$ and  \cite[Theorem (4.3.2)]{M} yields
$\Delta\nu_{8n-1}=\Delta\sigma_{16n-1}=0$ for $n\geq 1$.

Let $n\not\equiv 3\ (\bmod\ 4)$. By use of the exact sequence $(\mathcal{R}^{2n+1}_{2n+3})$
and \cite[p. 16, Table]{K}, we have  $\sharp\Delta\nu_{2n+1}=2$. So, by (\ref{JDel2})
and (\ref{W4}), the equality $(\ast)_{\nu_{2n+1}}$ holds if $n\neq 2^i-2$ with $i\geq 3$. By  $({\mathcal R}^{2n}_{2n+2})$, \cite[p. 161, Table]{K} and the fact that $\pi_{14}(SO(12))\cong\Z_{24}\oplus\Z_4$ \cite{ME},
$$
\sharp\Delta\nu_{2n}=\left\{ \begin{array}{ll} 24,
& \mbox{if}\ n\ \mbox{is even and $n\geq 4$};\\
12, & \mbox{if}\ n\ \mbox{is odd and $n\geq 3$ or $n=6$}.
\end{array}
\right.
$$
Hence, by (\ref{W4}), the equality $(\ast)_{\nu_n}$ for  $n\geq 4$ does not hold if and only if $n=2^i-3$ with $i\geq 4$.

By the fact that $\pi_k(SO(5))=0$ for $k=8,9$, the map $\Delta: \pi_k(\mathbb{S}^4)\to\pi_{k-1}(SO(4))\cong\pi_{k-1}(SO(3))\oplus\pi_{k-1}(\mathbb{S}^3)\cong(\Z_2)^2$ is an isomorphism.
This implies
$$
\sharp\Delta(\nu_4\eta_7)=\sharp\Delta((E\nu')\eta_7)=\sharp\Delta(\nu_4\eta_7^2)=\sharp\Delta((E\nu')\eta^2_7)=2.
$$
\par We obtain $[\iota_4,\nu_4\eta_7]=[\iota_4,(E\nu')\eta_7]=0$ and
$[\iota_4,(E\nu')\eta^2_7]=[\iota_4,\nu_4\eta^2_7]=0$ (\ref{hw}).

By $({\mathcal R}^6_{10})$ and the fact that $\pi_{10}(SO(6))\cong\Z_{120}\oplus\Z_2$, $\pi_{10}(SO(7))\cong\Z_8$ \cite[p. 162, Table]{K} and $\pi_{10}(\mathbb{S}^6)=0$,
we obtain $\sharp\Delta[\iota_6,\iota_6]=30$. On the other hand,
$\sharp[\iota_6,[\iota_6,\iota_6]]=3$.

We have
$\Delta\nu^2_{8n-1}=(\Delta\nu_{8n-1})\nu_{8n+1}=0$ for $n\geq 1$. By \cite[Lemma 3.1]{GM}, $\Delta\nu^2_{8n-k}=0$ if $n\geq 1$ and $k=3,4$.

Next, we need the formula parallel to (\ref{SO4}):
\begin{equation}\label{SO8}
\pm\Delta\iota_8=2[\iota_7]-[\eta_6]_8.
\end{equation}
Here, we use the fact that $\pi_7(SO(7))=\{[\eta_6]\}\cong\Z$ and
$\pi_7(SO(8))\cong\pi_7(\mathbb{S}^7)\oplus\pi_7(SO(7))\cong\Z^2$ generated by $[\iota_7]$ and $[\eta_6]_8$.
\par Let $n\equiv 3\ (\bmod\ 8)\geq 11$. We consider the commutative diagram
$$
(\diamond) \
\xymatrix{{}&{}&\pi_{n+6}(\S^n)\ar^{i'_*}[d]\ar^{\Delta}[dr]&\\
\pi_{n+6}(SO(n))
\ar[r]^(.5){{i_{n,n+3}}_\ast}&\pi_{n+6}(SO(n+3))
\ar[r]^{p_\ast}&
\pi_{n+6}(V_{n+3,3})\ar[d]^{p'_\ast}\ar[r]^(.5){\Delta'}&\pi_{n+5}(SO(n)
)\\
{}&{}&\pi_{n+6}(V_{n+3,2}),}
$$
where $V_{n,k}=SO(n)/SO(n-k)$ is the Stiefel manifold and the horizontal, vertical sequences are exact.

By \cite{BM},
${i_{n,n+3}}_\ast$ is a split epimorphism for $n\geq 19$ and in view of \cite[p. 16, Table]{K}, $\pi_{n+6}(SO(n+3))\cong\Z_2$.
For $n=11$, there exists an element $[\iota_7]_m\mu_8\in\pi_{17}(SO(m))$ for $m\geq 8$ such
that $J([\iota_7]_m\mu_8)=\sigma_m\mu_{m+7}$ and
$\sharp([\iota_7]_m\mu_8)=2$. This implies that
${i_{11,14}}_\ast:\pi_{17}(SO(11))\to\pi_{17}(SO(14))$ is a split epimorphism as well.

Denote by $\R P^n_k=\R P^n/\R P^{k-1}$ for $k\leq n$ the stunted real projective space.
We obtain  $\pi_{n+6}(V_{n+3,3})\cong\pi_{n+6}(\R P^{n+2}_n)\cong\pi^s_9(\R P^5_3)\cong\pi^s_7(\R P^2)\cong\Z_2$ and $\pi_{n+6}(V_{n+3,2})\cong\pi_{n+6}(\R P^{n+2}_{n+1})=0$.
Hence, by the diagram $(\diamond)$,
$\Delta\nu^2_n\neq 0$ and thus, by (\ref{W5}), we conclude that the equality $(\ast)_{\nu^2_n}$ for $n\ge 4$ does not hold if and only if $n=2^i-5$ with $i\geq 4$.

We recall the relation $\sigma'''=\{\nu_5,24\iota_8,\nu_8\}$.
For a lift $[\nu_5]\in\pi_8(SO(6))\cong\Z_{24}$ \cite{ME} of $\nu_5$, there exists a lift $[\sigma''']\in\{[\nu_5],24\iota_8,\nu_8\}\subset\pi_{12}(SO(6))\cong\Z_{60}$ of $\sigma'''$. This implies $\Delta\sigma'''=0$.
By (\ref{wsg''}), the equality $(\ast)_{\sigma''}$ holds.

In virtue of \cite[p.132, Table]{Mi}, $\pi_{14}(SO(7))\cong\Z_{2520}\oplus\Z_8\oplus\Z_2$ and $\pi_{14}(SO(9))\cong\Z_8\oplus\Z_2$.
We have $\pi_{14}(SO(8))\cong\pi_{14}(SO(7))\oplus\pi_{14}(\S^7)\cong\Z_{2520}\oplus\Z_8\oplus\Z_2\oplus\Z_{120}$.
By \cite[p. 21]{Ka}, $\Ker\{\Delta:\pi^8_{15}\to R^8_{14}\}=8\pi^8_{15}$. Therefore,
$$
\sharp\Delta(E\sigma')_{(2)}=\sharp\Delta(\sigma_8)_{(2)}=8.
$$
In virtue of \cite[Proposition IV.5]{Se} (\cite[(13.1)]{T}) and (\ref{JDel2}), we have
\begin{equation}\label{Serre}
\sharp\Delta\pi_{15}(\S^8;p)=p  \hspace{0.3cm} \mbox{for} \hspace{0.3cm} p=3, 5.
\end{equation}
This implies
\begin{equation}\label{DeEsig}
\sharp\Delta(E\sigma')=120 \ \mbox{and} \  \sharp\Delta\sigma_8=2520.
\end{equation}

On the other hand  $\sharp[\iota_8,E\sigma']=60$ (\ref{W60}) and $\sharp[\iota_8,\sigma_8]=120$ (\ref{W6}). Hence, the equality $(\ast)_{\alpha}$ for $\alpha=E\sigma', \sigma_8$ does not hold.



By \cite[p.\ 336]{KM}, $\sharp\Delta\sigma_{11}=2$ and we know $[\iota_{11},\sigma_{11}]=0$ (\ref{W6}). The remaining cases are for $\sigma_n$ with $n=9,10$ and $n\geq 12$. These are obtained by the parallel argument to $\nu_n$. This leads to (1).

By (\ref{JDel2}), the $J$-homomorphism restricted to $\{\Delta\alpha\}\subseteq\pi_{k-1}(SO(n))$ is a monomorphism provided the equality $(\ast)_\alpha$ holds.
Hence (1) implies (2). This completes the proof.
\end{pf}

Now, we show the following result extending Theorem \ref{PW}:
\begin{thm}\label{main}
The equality $G_{k+n}(\mathbb{R}P^n)={\gamma_n}_\ast G_{k+n}(\mathbb{S}^n)$ holds if $k\le 7$ except the following pairs: $(k,n)=(3,4),(4,4),(5,4),(6,4),(5,6),(7,8),(7,11),\\
(3,2^i-3)\; \mbox{with}\; i\ge 4\;\mbox{and}\;(6,2^i-5)\;\mbox{with $i\ge 5$}$.

Furthermore,

\mbox{\em (1)} $G_7(\mathbb{R}P^4)\supseteq 12\pi_7(\mathbb{R}P^4)$;

\mbox{\em (2)}
$G_{10}(\mathbb{R}P^4)\supseteq 3\pi_{10}(\mathbb{R}P^4)$;

\mbox{\em (3)}
$G_{11}(\mathbb{R}P^6)\supseteq 30\pi_{11}(\mathbb{R}P^6)$;

\mbox{\em (4)}
$G_{15}(\mathbb{R}P^8)\supseteq 2520\pi_{15}(\mathbb{R}P^8)$;

\mbox{\em (5)}
$G_{18}(\mathbb{R}P^{11})\supseteq 2\pi_{18}(\mathbb{R}P^{11})$;

\mbox{\em (5)}
$G_{2^i}(\mathbb{R}P^{2^i-3})\supseteq 2\pi_{2^i}(\mathbb{R}P^{2^i-3})$ for $i\ge 4$.
\end{thm}
\begin{pf}
In the light of Lemmas \ref{eg}, \ref{fund1}, \ref{nlem}, Proposition \ref{rc1}(1),
\cite{GM} and \cite{St}, we obtain the equality.

By (\ref{Dnu4}), $\sharp\Delta\nu_4=12$ and $\sharp\Delta(\nu^2_4)=3$.
This and (\ref{sl}) for $\mathbb{R}$ imply: (1)
$G_7(\mathbb{R}P^4)\supseteq
p'_\ast\pi_7(SO(5))=12\pi_7(\mathbb{R}P^4)$ and (2)
$G_{10}(\mathbb{R}P^4)\supseteq
p'_\ast\pi_{10}(SO(5))=3\pi_{10}(\mathbb{R}P^4)$. We know
$\sharp\Delta[\iota_6,\iota_6]=30$. Hence, Lemma \ref{fund1}(1)
implies (3).

\par By (\ref{DeEsig}),
$$
\Ker\{\Delta:\pi_{15}(\S^8)\to \pi_{14}(SO(8))\}=2520\pi_{15}(\S^8).
$$
This leads to (4).

(5) follows from the fact that
$\Ker\ \{\Delta:\pi_{18}(\mathbb{S}^{11})\to\pi_{17}(SO(11))\}=2\pi_{18}(\mathbb{S}^{11})$. (6) follows from the fact that
$\Ker\ \{\Delta:\pi_{2^i}(\mathbb{S}^{2^i-3})\to\pi_{2^i-1}(SO(2^i-3))\}=2\pi_{2^i}(\mathbb{S}^{2^i-3})$.
This completes the proof.
\end{pf}

In virtue of Proposition \ref{Gotc}, an upper bound of $G_7(\R P^4)$ is given as
${\gamma_4}_*G_7(\S^4)={\gamma_4}_*\{3[\iota_4,\iota_4],2E\nu'\}$ \cite[p. 7]{GM}. Hence,
\begin{rem}\label{upb} {\em
$6\gamma_4\nu_4\not\in G_7(\R P^4)$ and $\gamma_4(E\nu')\not\in G_7(\R P^4)$.}
\end{rem}

To state the next result, we recall from \cite{T} and \cite{GM}:
\begin{equation}\label{j1}
\Delta\varepsilon_{4n+3}=\Delta\mu_{4n+3}=0,
\end{equation}
\begin{equation}\label{j2}
\sharp[\iota_{8n},\varepsilon_{8n}]=\sharp[\iota_{8n},\eta_{8n}\sigma_{8n+1}]=\sharp[\iota_{8n},\bar{\nu}_{8n}]=\sharp [\iota_{8n},\mu_{8n}]=2,
\end{equation}
\begin{equation}\label{j3}
\sharp[\iota_{8n},\eta_{8n}\varepsilon_{8n+1}]=\sharp[\iota_{8n},\eta^2_{8n}\sigma_{8n+2}]=\sharp[\iota_{8n},\nu^3_{8n}]=2,
\end{equation}
\begin{equation}\label{j4}
[\iota_{8n},\alpha_{8n}]\ne[\iota_{8n},\mu_{8n}]\ \mbox{for}\ \alpha=\eta\varepsilon, \eta^2\sigma, \nu^3, 
\end{equation}
\begin{equation}\label{j5}
\sharp[\iota_{8n},\eta_{8n}\mu_{8n+1}]=2,\
\sharp[\iota_{8n},\beta_1(8n)]=3.
\end{equation}
Hence, by (\ref{W1})--(\ref{W6}), (\ref{j1})--(\ref{j5}), Proposition \ref{Gotc} and Corollary \ref{GRP}, we obtain:

\begin{prop}\label{8910}

\mbox{\em (1)}
$G_{n+k}(\mathbb{R}P^n)=\pi_{n+k}(\mathbb{R}P^n)$ if $n\equiv 3\ (\bmod\ 4)$ and $8\le k\le 10$.

\mbox{\em (2)}
$G_{n+k}(\mathbb{R}P^n)=0$ if $n\equiv 0\ (\bmod\ 8), \ n\ge 16$ and $8\le k\le 10$.
\end{prop}

\section{Gottlieb groups of complex and quaternionic projective spaces}
By \cite[Theorem 1]{G2}, we obtain
$$
G_2(\mathbb{C}P^n)=0 \hspace{5mm} \mbox{for} \hspace{5mm} n\ge 1
$$
proved in \cite{GG}, \cite{L} and \cite{LMW} as well.

We recall the group structures of $\pi_{2n-1+k}(SU(n))$ for $0\le k\le 4, 6$ and $k=7$ from \cite{BH}, \cite{Bo}, \cite{K}, \cite{Ma}, \cite{Ma1} and \cite{ME}.
$$
\pi_{2n-1}(SU(n))\cong\Z; \
\pi_{2n}(SU(n))\cong\Z_{n!} \hspace{5mm} \mbox{for} \hspace{5mm} n\geq 2;
$$
$$
\pi_{2n+1}(SU(n))\cong\left\{ \begin{array}{ll} \Z_2,
& \mbox{if}\;\;n\;\;\mbox{is even},\\
0, & \mbox{if}\;\;n\;\;\mbox{is odd};
\end{array}
\right.
$$
$$
\pi_{2n+2}(SU(n))\cong\left\{ \begin{array}{ll} \Z_{(n+1)!}\oplus\Z_2,
& \mbox{if}\;\;n\;\;\mbox{is even and}\ n\geq 4,\\
\Z_{(n+1)!/2}, & \mbox{if}\;\;n\;\;\mbox{is odd};
\end{array}
\right.
$$
$$
\pi_{2n+3}(SU(n))\cong\left\{ \begin{array}{ll} \Z_{(24,n)},
& \mbox{if}\;\;n\;\;\mbox{is even},\\
\Z_{(24,n+3)/2}, & \mbox{if}\;\;n\;\;\mbox{is odd and}\ n\geq 3;
\end{array}
\right.
$$
$$
\pi_{2n+5}(SU(n))\cong\left\{ \begin{array}{ll} \Z_{(24,n+1)},
& \mbox{if}\;\;n\;\;\mbox{is odd},\\
\Z_{(24,n+4)/2}, & \mbox{if}\;\;n\;\;\mbox{is even}; \end{array}
\right.
$$
$$
\pi_{2n+6}(SU(n))\cong\left\{ \begin{array}{ll} \pi_{2n+6}(SU(n+1)),
& \mbox{if}\;\;n\equiv 2,3\ (\bmod\ 4)\geq 3,\\
\pi_{2n+6}(SU(n+1))\oplus\Z_2, & \mbox{if}\;\;n\equiv 0,1\ (\bmod\ 4)\geq 4.
\end{array}
\right.
$$

Notice that $\pi_{2n}(SU(n))$ is generated by $\Delta_\C(\iota_{2n+1})=\omega_n(\C)=j_nE\gamma_{n-1}$, where
$j_n: E\mathbb{C}P^{n-1}\hookrightarrow SU(n)$ is the canonical inclusion satisfying $p_n j_n=Eq_{n-1}$. Then, by (\ref{qgam}),
$$
p_n \omega_n(\C)=(n-1)\eta_{2n-1}.
$$
Notice that $\pi_{4n+1}(SU(2n))$ is generated by $\Delta_\C(\eta_{4n+1})=\omega_{2n}\eta_{4n}$. By making use of the exact sequences $({\mathcal{C}}^n_{2n+k})$ for $0\leq k\leq 3, k=6$ and the results above, we obtain the following result which overlaps with that of \cite[pp. 163-5]{K}:
\begin{lem} \label{cc}
\mbox{\em (1)}
$\sharp\omega_n(\C)=n!$;

\mbox{\em (2)}
$\sharp\Delta_\C(\eta_{2n+1})=\left\{ \begin{array}{ll} 2,
& \mbox{if}\;\;n\;\;\mbox{is even},\\
1, & \mbox{if}\;\;n\;\;\mbox{is odd};
\end{array}
\right.
$

\mbox{\em (3)}
$\sharp\Delta_\C(\eta^2_{2n+1})=\left\{ \begin{array}{ll} 2,
& \mbox{if}\;\;n\;\;\mbox{is even},\\
1, & \mbox{if}\;\;n\;\;\mbox{is odd};
\end{array}
\right.
$

\mbox{\em (4)}
$\sharp\Delta_\C(\nu_{2n+1})=\left\{ \begin{array}{ll} (24,n), & \mbox{if}\;\;n\;\;\mbox{is even},\\
(24,n+3)/2, & \mbox{if}\;\;n\;\;\mbox{is odd and}\ n\geq 3;
\end{array}
\right.
$

\mbox{\em (5)}
$\sharp\Delta_\C(\nu^2_{2n+1})=\left\{ \begin{array}{ll} 1,
& \mbox{if}\;\;n\equiv 2,3\ (\bmod\ 4)\geq 2
,\\
2, & \mbox{if}\;\;n\equiv 0,1\ (\bmod\ 4)\geq 4.
\end{array}
\right.
$
\end{lem}

By Lemmas \ref{fund1}(1) and \ref{cc}(1), we obtain \cite[Theorem III.8]{L}:

\begin{thm}{\bf (Lang)}\hspace{2mm}\label{c4}
$n!\pi_{2n+1}(\mathbb{C}\P^n)\subseteq G_{2n+1}(\mathbb{C}P^n)$.
\end{thm}

By Theorems \ref{ex2}(1) and \ref{c4}, it holds

$$G_5(\mathbb{C}P^2)=P_5(\mathbb{C}P^2)=2\pi_5(\mathbb{C}P^2).$$

\bigskip

Now, we show:
\begin{thm}\label{GCP}
\mbox{\em (1)}\ Let $k=1,2$. Then:

$G_{k+2n+1}(\mathbb{C}P^n)
=\left\{
\begin{array}{ll}
0, & \mbox{if}\;\;n\;\;\mbox{is even},\\
\pi_{k+2n+1}(\mathbb{C}P^n)\cong\mathbb{Z}_2,&\; \mbox{if}\;\;n\;\;\mbox{is odd}.
\end{array}
\right.
$

\mbox{\em (2)}\\
$G_{2n+4}(\mathbb{C}P^n)\supseteq\left\{ \begin{array}{ll}
(24,n)\pi_{2n+4}(\mathbb{C}P^n)\cong\mathbb{Z}_{\frac{24}{(24,n)}}, & \mbox{if}\;\;n\;\;\mbox{is even},\\
\frac{(24,n+3)}{2}\pi_{2n+4}(\mathbb{C}P^n)\cong\mathbb{Z}_{\frac{48}{(24,n+3)}},
& \mbox{if}\;\;n\;\;\mbox{is odd}.
\end{array}
\right.
$

\noindent
In particular,  $G_{2n+4}(\mathbb{C}P^n)=2\pi_{2n+4}(\mathbb{C}P^n)$ \ if \ $n\equiv 2,10\ (\bmod\ 12)\geq 10$
except $n=2^{i-1}-2$ or $n\equiv 1,17\,(\bmod\,24)\ge 17$ and $G_{2n+4}(\mathbb{C}P^n)=\pi_{2n+4}(\mathbb{C}P^n)$ \ if \ $n\equiv 7,11\ (\bmod\ 12)$.

\mbox{\em (3)}
$G_{2n+6}(\mathbb{C}P^n)
=\pi_{2n+6}(\mathbb{C}P^n)
\cong\left\{
\begin{array}{ll}
0, & \mbox{if}\; n\ge 3,\\
\mathbb{Z}_2,& \mbox{if}\; n=2.
\end{array}
\right.
$

\mbox{\em (4)}
$G_{2n+7}(\mathbb{C}P^n)=\pi_{2n+7}
(\mathbb{C}P^n)$ if $n\equiv 2,3\ (\bmod\ 4)$.
\end{thm}

\begin{pf} By Proposition \ref{rc1}(1) and the fact that $G_{k+2n+1}(\S^{2n+1})=0$ if $k=1,2$ and $n$ is even, we have $G_{k+2n+1}(\mathbb{C}P^n)=0$ in this case.
By
Lemmas \ref{fund1}(1) and \ref{cc}(2)-(3), $G_{k+2n+1}(\mathbb{C}P^n)=\pi_{k+2n+1}(\mathbb{C}P^n)$ for $n$ odd. This leads to (1).

Notice that $G_6(\S^2)=\pi_6(\S^2)$ \cite{LMW}, \cite{GM}.
Lemmas \ref{fund1}(1)  and \ref{cc}(4) imply the first half of (2).

As it is easily seen (\ref{W4}), for $n$ even,
$(24,n)=\sharp[\iota_{2n+1},\nu_{2n+1}]=2$ if and only if $n\equiv
2,10\ (\bmod\ 12)\geq 10$ and $n\neq 2^{i-1}-2$. For $n$ odd,
$\frac{(24,n+3)}{2}=\sharp[\iota_{2n+1},\nu_{2n+1}]=2$
if and only if $n\equiv 1,17\ (\bmod\ 24)\geq 17$. Moreover,
$\frac{(24,n+3)}{2}=\sharp[\iota_{2n+1},\nu_{2n+1}]=1$ if and only if
$n\equiv 7,11\ (\bmod\ 12)$. This and Lemma \ref{fund1}(2) lead to the second half of (2).

Since $\Delta_\C: \pi_{10}(\S^5)\to\pi_9(SU(2))\cong\Z_3$ is trivial, we have $G_{10}(\C P^2)=\pi_{10}(\C P^2)\cong\Z_2$. This leads to (3).

(4) follows
from Lemmas \ref{fund1}(1) and \ref{cc}(5), and the proof is complete.
\end{pf}

We obtain:

\begin{cor} \label{GCP1}
\mbox{\em(1)}\
$G_k(\mathbb{C}P^2)=\pi_k(\mathbb{C}P^2)$ for $10\leq k\leq 12$;

\mbox{\em(2)}\ $G_k(\C P^2;p)=\pi_k(\C P^2;p)$ for an odd prime $p$;

\mbox{\em (3)}\
$G_k(\C P^{2n+1})\supseteq\gamma_{2n+1}\eta^m_{4n+3}\circ\pi_k(\S^{4n+3+m})$ for $m=1,2$;

\mbox{\em (4)}\
\mbox{\em (i)}
$G_k(\C P^{2n})\supseteq 2(12,n)\gamma_{2n}\nu_{4n+1}\circ\pi_k(\S^{4n+4})$;

\hspace{.4cm}
\mbox{\em (ii)}
$G_k(\C P^{2n+1})\supseteq (12,n+2)\gamma_{2n+1}\nu_{4n+3}\circ\pi_k(\S^{4n+6})$;

\mbox{\em (5)}\
$G_k(\C P^{4n+2})\supseteq\gamma_{4n+2}\nu^2_{8n+5}\circ\pi_k(\S^{8n+11})$;

\mbox{\em (6)}\
$G_{8n+11}(\C P^{4n+1})=\pi_{8n+11}(\C P^{4n+1})$\ if\ $n\geq 2$;

\mbox{\em (7)}\
$G_{8n+k}(\C P^{4n+3})=\pi_{8n+k}(\C P^{4n+3})$ for  $k=28,29$\ if\ $n\geq 2$.
\end{cor}
\begin{pf}
By
\cite[Lemma 3.1.i)]{MT1}, we get $\Delta_\mathbb{C}(\pi_k(\mathbb{S}^5))=0$ for $10\leq k\leq 12$.
Hence, Lemma \ref{fund1}(1) yields (1).

By Lemma \ref{eg} and Theorem \ref{GCP}, we have (2)-(5).
By Lemma \ref{cc}(2), $\Delta_\C(\eta_{4n+3})=0$. Since $\varepsilon_{4n+3}\in\{\eta_{4n+3},\nu_{4n+4},2\nu_{4n+7}\}\ (\bmod\ \eta_{4n+3}\sigma_{4n+4})$, we have
$$
\Delta_\C\varepsilon_{4n+3}\in\{\omega_{2n+1},\eta_{4n+2},\nu_{4n+3}\}\circ 2\nu_{4n+7}\subseteq\pi_{4n+7}(SU(2n+1))\circ 2\nu_{4n+7}.
$$
Then, by the group structure $\pi_{4n+7}(SU(2n+1))\cong\Z_{2(12,n+1)}$, we obtain
$\Delta_\C\varepsilon_{4n+3}=0$ for $n$ even. This leads to (6).

By Lemma \ref{cc}(4), $\sharp\Delta_\C(\nu_{8n+7})=(12,2n+3)$. We recall from \cite{Mi0} that $\pi^s_{21}(\S^0)=\{\eta\bar{\kappa},\nu\nu^\ast\}\cong(\Z_2)^2$ and $\pi^s_{22}(\S^0)=\{\eta^2\bar{\kappa},\nu\bar{\sigma}\}\cong(\Z_2)^2$. This implies $\Delta_\C(\eta_{4n+3}\bar{\kappa}_{4n+4})=\Delta_\C(\nu_{8n+7}\nu^\ast_{8n+10})=0$ and leads to (7) for $k=28$. By the parallel argument, we have the case $k=29$.
\end{pf}

Next, we consider the case of the quaternionic projective space.
By \cite[Theorem 1]{G2}, we obtain
$$
G_4(\mathbb{H}P^n)=0 \hspace{5mm} \mbox{for} \hspace{5mm} n\ge 1
$$
proved in \cite{GG}, \cite{L} and \cite{LMW} as well.

By Proposition \ref{exh1}(1);(4) and \cite{GM},
$$
G_k(\mathbb{H}P^n)=0 \hspace{3mm} \mbox{if} \ n\geq 1 \ \mbox{and} \ k=5,6;
$$
$$
G_k(\mathbb{H}P^n)=0 \hspace{3mm} \mbox{if}\ n\geq 3 \ \mbox{and} \ k=12,13;
$$

We recall the group structures $\pi_{4n-1+k}(Sp(n))$ for $k=-1,0$ and $2\leq k\leq 6$
from \cite{Bo}, \cite{MT2}, \cite{Mi1} and \cite{ME}:
$$
\pi_{4n-2}(Sp(n))=0; \ \pi_{4n-1}(Sp(n))\cong\Z;
$$
$$
\pi_{4n+1}(Sp(n))\cong\left\{\begin{array}{ll} 0,& \mbox{if}\;\;n\;\;\mbox{is even},\\
\Z_2, & \mbox{if}\;\;n\;\;\mbox{is odd};
\end{array}
\right.
$$
$$
\pi_{4n+2}(Sp(n))\cong\left\{\begin{array}{ll} \Z_{(2n+1)!},& \mbox{if}\;\;n\;\;\mbox{is even},\\
\Z_{2(2n+1)!}, & \mbox{if}\;\;n\;\;\mbox{is odd};
\end{array}
\right.
$$
$$
\pi_{4n+3}(Sp(n))\cong\Z_2;
$$
$$
\pi_{4n+4}(Sp(n))\cong\left\{\begin{array}{ll} (\Z_2)^2,& \mbox{if}\;\;n\;\;\mbox{is even},\\ \Z_2,& \mbox{if}\;\;n\;\;\mbox{is odd};\end{array}\right.
$$
$$
\pi_{4n+5}(Sp(n))\cong\left\{\begin{array}{ll} \Z_{(24,n+2)}\oplus\Z_2,& \mbox{if}\;\;n\;\;\mbox{is even};\\
\Z_{(24,n+2)}, & \mbox{if}\;\;n\;\;\mbox{is odd};
\end{array}
\right.
$$

By making use of the exact sequences $(\mathcal{\mathcal{SH}}^n_{4n+k})$ for $2\leq k\leq 5$ and the results above, we obtain:
\begin{lem} \label{hh}
\mbox{\em (1)}
$\sharp\Delta_\H(\iota_{4n+3})=\left\{\begin{array}{ll} (2n+1)!,& \mbox{if}\;\;n\;\;\mbox{is even},\\
2(2n+1)!, & \mbox{if}\;\;n\;\;\mbox{is odd};
\end{array}
\right.
$

\mbox{\em (2)}
$\sharp\Delta_\H(\eta^k_{4n+3})=2$ for $k=1,2$;

\mbox{\em (3)}
$\sharp\Delta_\H(\nu_{4n+3})=(24,n+2).$
\end{lem}

Notice that $\pi_{4n+2}(Sp(n))$ is generated by $\Delta_\H(\iota_{4n+3})=\omega_n(\H)$.

By Lemmas \ref{fund1}(1), \ref{hh}(1);(3) and \ref{eg}, we obtain:
\begin{thm}\label{hthm}
\mbox{\em (1)}
$G_{4n+3}(\mathbb{H}P^n)\supseteq\left\{\begin{array}{l}(2n+1)!{\gamma_n}_\ast\pi_{4n+3}(\S^{4n+3}),\,if\, n\,is\, even,\\
2(2n+1)!{\gamma_n}_\ast\pi_{4n+3}(\S^{4n+3}),\, if\,n\,is\, odd;
\end{array}
\right.
$

 \mbox{\em (2)}
$G_k(\H P^n)\supseteq(24,n+2)\gamma_n\nu_{4n+3}
\circ\pi_k(\S^{4n+6})$. In particular, we derive $G_{4n+6}(\mathbb{H}P^n)\supseteq(24,n+2){\gamma_n}_\ast\pi_{4n+6}(\S^{4n+3})\cong\Z_{\frac{24}{(24,n+2)}}$ for $n\geq 2$;
\end{thm}

Theorem \ref{hthm}(1) improves \cite[Corollary 2.7]{LMW}. We obtain:

\begin{cor}\label{GHP1}
{\em (1)}
$G'_{8n+9}(\mathbb{H}P^{2n})=0$,\ $G'_{8n+10}(\mathbb{H}P^{2n+1})={\gamma_{2n+1}}_\ast\pi_{8n+10}(\S^{8n+7})$\ for $n\not\equiv 0\ (\bmod\ 3)$
 \ and \ $G'_{8n+13}(\mathbb{H}P^{2n+1})
={\gamma_{2n+1}}_\ast\pi_{8n+13}(\S^{8n+7})$.

{\em (2)}
$G'_{4n+14}(\mathbb{H}P^n)={\gamma_n}_\ast\pi_{4n+14}(S^{4n+3})$\ for $n\equiv 5,9\ (\bmod\ 12), n\ge 5$ and
$n\equiv 15,23\ (\bmod\ 24)$.
\end{cor}
\begin{pf}
(1) is a direct consequence of Theorem \ref{hthm}(2); (3) and Theorem \ref{PHP1}(4).

We know  $\zeta_n\in\{\nu_n,\sigma_{n+3},16\iota_{n+10}\}$ for $n\ge 13$. Moreover, $\pi_{4n+13}(Sp(n))\cong\mathbb{Z}_8$ if $n\equiv 1\ (\bmod\ 4)$ and $n\ge 5$ and $\cong\mathbb{Z}_{(128,4(n-3))}$ if $n\equiv 3\ (\bmod\ 4)$ (\cite{Mo1}, \cite{Mo2}). So, for $n$ satisfying $(24,n+2)=1, n\equiv 1\ (\bmod\ 4)$ or $n\equiv 7\ (\bmod\ 8)$, we obtain
$$
\Delta_\mathbb{H}\zeta_{4n+3}\in\{\omega_n(\H),\nu_{4n+2},\sigma_{4n+5}\}\circ 16\iota_{4n+13}
$$
$$
\subseteq\pi_{4n+13}(Sp(n))\circ16\iota_{4n+13}=0.
$$
This and Lemma \ref{fund1} lead to (2).
\end{pf}

By use of the exact sequence $(\mathcal{\mathcal{SH}}^2_{11})$, Lemma \ref{hh}(2)  and the fact that $\pi_{12}(Sp(2))=\{i\varepsilon_3\}\cong\Z_2$\ for the inclusion\ $i:\S^3\hookrightarrow Sp(2)$ \cite[Theorem 5.1]{MT1}, we obtain
$$
\Delta_\H(\eta_{11})=i\varepsilon_3.
$$
So, by the relation $\varepsilon_3\sigma_{11}=0$ \cite[Lemma 10.7]{T}, we obtain:
\begin{exam}\label{exa} {\em
$G_{18+k}(\mathbb{H}P^2)\supseteq
\{\gamma_2\eta^k_{11}\sigma_{11+k}\}\cong\Z_2$\ for $k=1,2$.}
\end{exam}

We have $\nu\rho\in\nu\circ\langle\sigma,2\sigma,8\iota\rangle=8\langle\nu,\sigma,2\sigma\rangle=8\nu^\ast=0$. By \cite[(10.22)]{T}, $H(\nu_{10}\rho_{13})=E(\nu_9\wedge\nu_9)\circ 4\nu_{25}=0$. So, by the fact that $\nu_{10}\rho_{13}\in E\pi^9_{27}=\{\bar{\zeta}_{10},\bar{\sigma}_{10}\}\cong \pi^s_{19}$, we obtain
\begin{equation}\label{nurh}
\nu_{10}\rho_{13}=0.
\end{equation}
By \cite[Lemma 6.1]{Mi0} and its proof, we have
$$
\eta_6\bar{\kappa}_7\in\{\nu^2_6,2\iota_{12},\kappa_{12}\}\ (\bmod\ \nu^2_6\circ\pi^{12}_{27}+\pi^6_{13}\circ\kappa_{13})
$$
$$
=\nu^2_6\circ\{E^3\rho',\bar{\varepsilon}_{12}\}+\{\sigma''\kappa_{13}\}.
$$
We obtain $\nu^2_6\bar{\varepsilon}_{12}=\nu^2_6\eta_{12}\kappa_{13}=0$. By \cite[Lemma 6]{Os}, $\nu^2_6(E^3\rho')=0$.
We know $\sigma''\kappa_{13}\in\pi^6_{27}=\{\eta_6\bar{\kappa}_7\}$, $E(\sigma''\kappa_{13})=2\sigma'\kappa_{14}=0$ and $\eta\bar{\kappa}\ne 0$ \cite[Theorem A]{Mi0}. Hence
\begin{equation}\label{MOs}
\eta_6\bar{\kappa}_7=\{\nu^2_6,2\iota_{12},\kappa_{12}\}.
\end{equation}

As it is well known \cite[(2.10)(b)]{J3}, for the projection $p_n=p_{n,\H}: Sp(n)\to\S^{4n-1}$, we have
\begin{equation}\label{pom}
{p_n}_\ast\omega_n(\H)=\pm(n+1)\nu_{4n-1}\ (n\ge 2).
\end{equation}

Then, we obtain:

\begin{lem}\label{spnu}
$\Delta_\H(\eta_n\bar{\kappa}_{n+1})=0$
if\ $n\equiv 15\ (\bmod\ 16)$.
\end{lem}
\begin{pf}
We use $\omega_n=\omega_n(\H)$.
By (\ref{nurh}) and (\ref{MOs}),
it holds that
$$
\{\nu^2_n,2\iota_{n+6},\kappa_{n+6}\}\ni\eta_n\bar{\kappa}_{n+1}\ (\bmod\ \sigma_n\kappa_{n+7})\ \mbox{for}\ n\geq 10.
$$
Set $n=4m$. By Lemma \ref{hh}(3),
$\sharp(\omega_n\nu_{4n+2})=2(12,2m+1)$. So we can define the Toda bracket
$$
\{\omega_n\nu_{4n+2},2(12,2m+1)\iota_{4n+5},\kappa_{4n+5}\}\subset\pi_{4n+20}(Sp(n)).
$$
By (\ref{pom}), we obtain
\begin{eqnarray*}
{p_n}_*\{\omega_n\nu_{4n+2},2(12,2n+1)\iota_{4n+5},\kappa_{4n+5}\}
&\subseteq&\{\nu^2_{4n-1},2(12,2n+1)\iota_{4n+5},\kappa_{4n+5}\}\\
&\ni&\eta_{4n-1}\bar{\kappa}_{4n}\ (\bmod\
\sigma_{4n-1}\kappa_{4n+6}).
\end{eqnarray*}
So, by the relation $\sigma_{4n-1}\kappa_{4n+6}=0$ for $n\ge 4$ \cite[Proposition 7.2]{Mi0}, we obtain
$$
{p_n}_*\{\omega_n\nu_{4n+2},2(12,2m+1)\iota_{4n+5},\kappa_{4n+5}\}=\eta_{4n-1}\bar{\kappa}_{4n}.
$$
This completes the proof.
\end{pf}

By \cite[Example 4]{T}, we know $\sigma^2\in\langle\nu,\sigma,\nu\rangle\
(\bmod\ \kappa)$. Since $\{\eta_{10},\nu_{11},\sigma_{14}\}=2[\iota_{10},\nu_{10}]$ \cite{Mu2} and
$\eta_{10}\kappa_{11}\ne 0$, we obtain
\begin{equation}\label{tnusg11}
\{\nu_{11},\sigma_{14},\nu_{21}\}=a\sigma^2_{11}\ \mbox{for $a$ odd}.
\end{equation}
By (\ref{tnusg11}) and the definition of $\xi_{12}$, we have
$$
\sigma^3_{15}\in\{\nu_{15},\sigma_{18},\nu_{25}\}\circ\sigma_{29}=-\nu_{15}\circ\{\sigma_{18},\nu_{25},\sigma_{28}\}
$$
$$
\ni-\nu_{15}\xi_{18}\ (\bmod\ \nu_{15}\circ\pi^{18}_{29}\circ\sigma_{29}=0).
$$
Hence, we have
\begin{equation}\label{sg315}
\sigma^3_{15}=\nu_{15}\xi_{18}.
\end{equation}

From the fact that $\sharp\Delta_\H(\nu_{8n+7})=(24,2n+3)$ and (\ref{tnusg11}), we have
$$
\Delta_\H(\sigma^2_{8n+7})\in\{\omega_{2n+1},(24,2n+3)\nu_{8n+6},\sigma_{8n+9}\}\circ\nu_{8n+17}
$$
and hence,
\begin{equation}\label{Dsgm3}
\Delta_\H(\sigma^3_{8n+7})=0.
\end{equation}

We show:
\begin{cor} \label{HGP2}
$G'_{16n+k}(\H P^{4n-1})={\gamma_{4n-1}}_\ast\pi_{16n+k}(\S^{16n-1})$ for $k=20; n\ge 1$ and $k=21; n\ge 2$.
\end{cor}
\begin{pf}
By \cite[(7.10)]{Mi0}, $E^3(\lambda\nu_{31})=[\iota_{16},\nu^2_{16}]$. So, by (\ref{sg315}), \cite[Corollary 12.25]{T} and (\ref{W6}),
$$
\sigma^3_{16}=\nu_{16}\xi_{19}=\nu_{16}\circ([\iota_{19},\iota_{19}]-\nu^*_{19})=E^3(\lambda\nu_{31})-\nu_{16}\nu^*_{19}
$$
and $\sigma^3_{15}+\nu_{15}\nu^*_{18}-E^2(\lambda\nu_{31})\in P\pi^{31}_{38}=0.$

By Lemma \ref{hh}(3) and (\ref{Dsgm3}), we have $\Delta_\H((E^2\lambda)\nu_{33})=\Delta_\H(\nu_{15}\nu^*_{18})+\Delta_\H(\sigma^3_{15})=0$ and
$G_{16n+20}(\H P^{4n-1})\ni{\gamma_{4n-1}}_\ast\sigma^3_{16n-1}$. By Lemma \ref{spnu},
$G_{16n+20}(\H P^{4n-1})\ni{\gamma_{4n-1}}_\ast\eta_{16n-1}\bar{\kappa}_{16n}$. This leads to the assertion for $k=20$.

By Theorem \ref{hthm}(3), $G_{8n+21}(\H P^{2n-1})\ni{\gamma_{2n-1}}_\ast\nu_{8n-1}\bar{\sigma}_{8n+2}$. By Lemma \ref{spnu},
$$
\Delta_\H(\eta^2_{16n-1}\bar{\kappa}_{16n+1})=\Delta_\H(\eta_{16n-1}\bar{\kappa}_{16n})\circ\eta_{16n+20}=0.
$$
This and the group structure of $\pi_{16n+21}(\S^{16n-1})$ for $n\ge 2$ \cite[Theorem B]{Mi0} lead to the assertion for $k=21$. This completes the proof.
\end{pf}

Finally, we show:
\begin{prop}\label{CHP2ex2} {\em (1)} $G_{18}(\mathbb{H}P^2)\supseteq 40{\gamma_2}_\ast\pi_{18}(\S^{11})$;

{\em(2)} $G_{21}(\mathbb{H}P^2)\supseteq 2{\gamma_2}_\ast\pi_{21}(\S^{11})$;

{\em(3)} $G_{22}(\mathbb{H}P^2)\supseteq 8{\gamma_2}_\ast\pi_{22}(\S^{11})\cong\Z_{63}$;

{\em(4)} $G_{22}(\mathbb{H}P^3)\supseteq 4{\gamma_3}_\ast\pi_{22}(\S^{15})\cong\Z_{60}$. \end{prop}
\begin{pf} By \cite[Theorem 5.1]{MT1} and \cite[Lemmas 3.1-4]{MT2}, we see that

(1) $\Ker\{\Delta_{\H}: \pi_{18}(\S^{11})\to\pi_{17}(Sp(2))\}=40\pi_{18}(\S^{11})$;

(2) $\Ker\{\Delta_{\H}:\pi_{21}(\S^{11})\to\pi_{20}(Sp(2))\}=2\pi_{21}(\S^{11})$;

(3) $\Ker\{\Delta_{\H}: \pi_{22}(\S^{11})\to\pi_{21}(Sp(2))\}=8\pi_{22}(\S^{11})$;

(4) $\Ker\{\Delta_{\H}:\pi_{22}(\S^{15})\to\pi_{21}(Sp(3))\}=4\pi_{22}(\S^{15})$.

\noindent Hence, by Lemma \ref{fund1}(1), we obtain the assertion.  This completes the proof.
\end{pf}

\noindent
{\bf Problem.}\ Find a pair $(k,n)$ satisfying $G_k(\H P^n)\cap{i_\H}_*E\pi_{k-1}(\S^3)\ne 0$.

\section{The case of the Cayley projective plane} Let $\mathbb{K}P^2=\S^8\cup_{\sigma_8}e^{16}$ be the Cayley projective plane. By \cite[Theorem 1]{G2}, we obtain
$$
G_8(\mathbb{K}P^2)=0
$$
proved in \cite{GG} and \cite{LMW} as well.

We recall the fibration
$$
p : F_4\longrightarrow F_4/Spin(9)=\mathbb{K}P^2
$$
and the induced exact sequence:
\begin{equation}\label{FSpine}
\cdots\rarrow{}
\pi_k(F_4)\rarrow{p_*}\pi_k(\K P^2)\rarrow{\Delta_\K}\pi_{k-1}(Spin(9))\rarrow{i_*}\cdots.
\end{equation}
By (\ref{FSpine}) and Lemma \ref{SL}, we have
\begin{equation}\label{FSpine2}
\Ker\ \Delta_\K\subseteq G_k(\K P^2).
\end{equation}

Let $i_{\mathbb{K}}: \mathbb{S}^8\hookrightarrow\mathbb{K}P^2$
be the inclusion map and $p_\K: \mathbb{K}P^2\to\S^{16}$ the collapsing map. We recall
$$
\pi_k(\S^8)={\sigma_8}_\ast\pi_k(\S^{15})\oplus E\pi_{k-1}(\S^7)
$$
and
the relation
\begin{equation}\label{sn}
 \sigma'\nu_{14}=x\nu_7\sigma_{10} \; \mbox{for some odd $x$ \; \cite[(7.19)]{T}}.
\end{equation}

First, we show:
\begin{lem}\label{JT}
$\pi_n(\mathbb{K}P^2)={i_{\mathbb{K}}}_\ast E\pi_{n-1}(\S^7)\cong\pi_{n-1}(\S^7)$ for $n\leq 21$ or $n=25,27,28$;

$\pi_{22}(\mathbb{K}P^2)=\{{i_{\mathbb{K}}}\kappa_8\}\cong\Z_4$ and $i_\K(E\sigma')\sigma_{15}=0$;

$\pi_{23}(\mathbb{K}P^2)\cong\Z\oplus\Z_{120}\{i_{\mathbb{K}}E\rho''\}\oplus\Z_2\{i_{\mathbb{K}}\bar{\varepsilon}_8\}
\oplus\Z_2\{i_{\mathbb{K}}(E\sigma')\varepsilon_{15}\}$ and
$i_{\mathbb{K}}((E\sigma')\varepsilon_{15})=i_{\mathbb{K}}((E\sigma')\bar{\nu}_{15})$;

$\pi_{24}(\mathbb{K}P^2)={i_{\mathbb{K}}}_*\{(E\sigma')\mu_{15},\mu_8\sigma_{15},\eta_8\bar{\varepsilon}_9\}\cong(\Z_2)^3$.

There is the short exact sequence
$$0\to \mathbb{Z}_{24}\oplus\mathbb{Z}_2\to\pi_{26}(\K P^2)\to\mathbb{Z}_{24}\to 0$$
and $\pi_{26}(\K P^2;2)\cong\Z_{64}\oplus\Z_2$ with a generator of $\Z_{64}$
given by $\tilde{\nu}_{21}\in \{i_\K,\sigma'\sigma_{14},\nu_{21}\}_1$.
\end{lem}
\begin{pf}
Although the result except for $n=24,25,27$ is obtained by \cite[Theorem 7.2]{Mi}, we give a proof. By the homotopy exact sequence of a pair $(\mathbb{K}P^2,\S^8)$ and the Blakers-Massey theorem, the assertion for $n\leq 21$ is obtained.

Let $n\leq 28$. Then, by \cite[Proposition 7.1]{Mi}, the $2$- and $3$-primary components of  $\pi_n(\mathbb{K}P^2)$ are obtained by determining $\pi_{n-1}(Y)$, where $Y$ is the $28$-skeleton of the loop space $\Omega(\K P^2)$:
$$
Y=\S^7\cup_{\sigma'\sigma_{14}+\alpha'}e^{22}\ \mbox{with $E\alpha'=[[\iota_8,\iota_8],\iota_8]$}.
$$
Recall  $\pi_{21}(\S^7)=\{\sigma'\sigma_{14},\kappa_8\}\cong\Z_{24}\oplus\Z_4$.
So, the group $\pi_{21}(Y)=\{i_\K\kappa_7\}\cong\Z_4$ is obtained the homotopy exact sequence of a pair $(Y,\S^7)$.
Hence, the group $\pi_{22}(\mathbb{K}P^2)$ and the relation
$i_\mathbb{K}(E\sigma')\sigma_{15}=0$ follows.

Next, recall
$\pi_{22}(\S^7)=\{\rho'',\sigma'\bar{\nu}_{14},\sigma'\varepsilon_{14},
\bar{\varepsilon}_7\}\cong\Z_{120}\oplus(\Z_2)^3$. By the Blakers-Massey
Theorem, $\pi_{22}(Y,\S^7)\cong\Z$ and $\pi_{23}(Y,\S^7)
\cong\pi_{23}(\S^{22})$. Since $(\sigma'\sigma_{14}+\alpha')\eta_{21}=
\sigma'(\bar{\nu}_{14}+\varepsilon_{14})$, we obtain the group
$\pi_{22}(Y)$ and the relation $i_{\mathbb{K}}((E\sigma')\varepsilon_{15})=i_{\mathbb{K}}((E\sigma')\bar{\nu}_{15})$.

By the parallel argument, we have the group $\pi_{23}(Y)$. This leads to the assertion for $n=23,24$.

We also recall
$\pi_{24}(\S^7)=\{(E\sigma')\eta_{15}\mu_{16},\nu_8\kappa_{11},
\bar{\mu}_8,\eta_8\mu_9\sigma_{18}\}\cong(\Z_2)^4$. By the Blakers-Massey
Theorem,  $\pi_{24}(Y,\S^7)\cong\pi_{24}(\S^{22})$.
In view of the relations $\sigma_{14}\eta^2_{21}=\eta_{14}\varepsilon_{15}+\nu^3_{14},
\sigma'\nu^3_{14}=\nu_7\sigma_{10}\nu^2_{17}=\eta_7\bar{\varepsilon}_8$\
(\ref{sn}), \cite[Lemma 12.10]{T},
we obtain $\sigma'\sigma_{14}\eta^2_{21}=E\zeta'+\eta_7\bar{\varepsilon}_8\neq 0$
\cite[Theorem 12.6, (12.4)]{T}.
Hence, $\partial: \pi_{24}(Y,\S^7)\to\pi_{23}(\S^7)$ is a monomorphism.
Since $\pi_{25}(Y,\S^7)\cong\pi_{25}(\S^{22})$ and $(\sigma'\sigma_{14}+\alpha')\nu_{21}=\alpha'\nu_{21}=0$,
we obtain $\partial\pi_{25}(Y,\S^7)=0$. This leads to $\pi_{25}(\K P^2)\cong\pi_{24}(\S^7)$.

Since $\pi_{25+k}(Y,\S^7)\cong\pi_{25+k}(\S^{22})=0$ for $k=1,2$, we get that $i_*: \pi_{26}(\S^7)\to\pi_{26}(Y)$ is an isomorphism.
This leads to the assertion for $n=27$ and completes the proof.
\par Finally, by $\pi_{27}(\S^{22})=0$ and $(\sigma'\sigma_{14}+\alpha')\nu^2_{21}=0$,
\cite[(7.6)]{Mi} yields $\pi_{28}(\K P^2)\cong\pi_{27}(\S^7)$.

Because $\pi_{26}(\S^{22})=0$, $\pi_{25}(\S^7)\cong\mathbb{Z}_{24}\oplus\mathbb{Z}_2$,
$\pi_{25}(\S^{22})\cong\mathbb{Z}_{24}$ and $(\sigma'\sigma_{14}+\alpha')\nu_{21}=0$
\cite[(7.6)]{Mi} lead to the short exact sequence
$$0\to \mathbb{Z}_{24}\oplus\mathbb{Z}_2\to\pi_{26}(\K P^2)\to\mathbb{Z}_{24}\to 0.$$
Observe $8\{i_\K, (E\sigma')\sigma_{15}, \nu_{22}\}=-i_\K \{(E\sigma')\sigma_{15},
\nu_{22}, 8\iota_{25}\} \supseteq -i_\K(E\sigma') \{\sigma_{15},
\nu_{22}, 8\iota_{25}\} \ni -i_\K(E\sigma')\zeta_{15}=
xi_\K\zeta_8\sigma_{19}\ (\bmod \ 8{i_\K}_*\pi^8_{26}=0)$
for some an odd integer $x$ (\cite[Lemma 12.12]{T}).
Consequently, we derive that $8\tilde{\nu}_{21}=xi_\K\zeta_8\sigma_{19}$ and the proof is complete.
\end{pf}

\bigskip

Since $\pi_7(Spin(9))\cong\mathbb{Z}$ and $\pi_7(F_4)=0$ \cite[p. 132, Table]{Mi}, (\ref{FSpine}) implies
\begin{equation}\label{iKb}
\Delta_\mathbb{K}(i_\mathbb{K})=b \ \mbox{for a generator of $\pi_7(Spin(9)$}.
\end{equation}

For any element $\beta\in\pi_m(\K P^2)$ with $m\ge 8$, we have
$$
[i_\K(E\alpha),\beta]=[i_\K,\beta]\circ E^m\alpha.
$$
In particular, for $\alpha\in\pi^7_{n-1}$ satisfying $E^8\alpha=0$, we obtain
\begin{equation}\label{Kcyc}
i_\K(E\alpha)\in P_n(\K P^2).
\end{equation}

Now, we show:
\begin{thm} \label{o2}
\mbox{\em(1)}
$P_k(\mathbb{K}P^2)=0$ for $k=9,10,12,13,14,20$;

\mbox{\em(2)}
$G_{11}(\mathbb{K}P^2)=P_{11}(\mathbb{K}P^2)=
8\pi_{11}(\mathbb{K}P^2)\cong\Z_3$;

\mbox{\em(3)}
$8\pi_{15}(\mathbb{K}P^2)\subseteq G_{15}(\mathbb{K}P^2)$;

\mbox{\em(4)}
${i_\K}_* E\{\sigma'\eta_{14}\}\subseteq P_{16}(\K P^2)\subseteq{i_\K}_* E\{\sigma'\eta_{14},\eta_7\sigma_8\}$;

\mbox{\em(5)}
${i_\K}_* E\{\sigma'\eta^2_{14}\}\subseteq P_{17}(\K P^2)\subseteq
{i_\K}_* E\{\sigma'\eta^2_{14},\nu^3_7+\eta_7\varepsilon_8\}={i_\K}_* E\{\sigma'\eta_{14},\eta_7\sigma_8\}\circ\eta_{15}$;

\mbox{\em(6)}
$8\pi_{18}(\mathbb{K}P^2)\subseteq G_{18}(\mathbb{K}P^2)\subseteq
P_{18}(\mathbb{K}P^2)={i_{\mathbb{K}}}_*\{\nu_8\sigma_{11},\beta_1(8)\}\cong\Z_{24}$;

\mbox{\em(7)}
${i_\K}_* E\{\bar{\nu}_7\nu_{15}\}\subseteq P_{19}(\K P^2)={i_\K}_* E\{\bar{\nu}_7\nu_{15}\}+8\pi_{19}(\K P^2)$;

\mbox{\em(8)}
$2\pi_{21}(\mathbb{K}P^2)\subseteq G_{21}(\mathbb{K}P^2)\subset P_{21}(\mathbb{K}P^2)=\pi_{21}(\mathbb{K}P^2)$.
\end{thm}
\begin{pf}
Since $\pi_k(\mathbb{K}P^2)\cong\pi_{k-1}(\S^7)=0$ for $k=12,13,20$ by Lemma \ref{JT}, the assertion for $k=12,13,20$ of (1) follows.

By Lemma \ref{JT},
$$
[i_{\mathbb{K}},i_\mathbb{K}\eta_8]=i_{\mathbb{K}}[\iota_8,\eta_8]=i_\mathbb{K}(E\sigma')\eta_{15}\neq 0,
$$
$$
[i_{\mathbb{K}},i_\mathbb{K}\eta^2_8]=i_\mathbb{K}(E\sigma')\eta^2_{15}\neq 0,
$$
$$
[i_{\mathbb{K}},i_\mathbb{K}\nu^2_8]=i_\mathbb{K}(E\sigma')\nu^2_{15}=i_\mathbb{K}\nu_8\sigma_{11}\nu_{18}\neq 0.
$$
This leads to the assertion for $k=9, 10, 14$ of (1).

We recall $\pi_{11}(\K P^2)=\{i_\K\nu_8\}\cong\Z_{24}$.
By \cite[Proposition 3.1]{T}, $\alpha_2(3)\wedge\alpha_1(3)=\alpha_2(6)\alpha_1(13)=-\alpha_1(6)\alpha_2(9).$
Hence, by (\ref{12beta}),
$$
\alpha_2(7)\alpha_1(14)=\alpha_1(7)\alpha_2(10)=0.
$$
By the relation $\pm[\iota_8,\nu_8]=2\sigma_8\nu_{15}-(E\sigma')\nu_{15}$ and Lemma \ref{JT}, $\sharp(i_\mathbb{K}[\iota_8,\nu_8])=\sharp(i_\mathbb{K}(E\sigma')\nu_{15})=8$ and $i_\mathbb{K}[\iota_8,\alpha_1(8)]=\pm i_\K\alpha_2(8)\alpha_1(15)=0$.
This implies $P_{11}(\mathbb{K}P^2)=0$ or $P_{11}(\mathbb{K}P^2)=\{i_\mathbb{K}\alpha_1(8)\}\cong\Z_3$.
On the other hand, by (\ref{iKb}) and the fact that
$\pi_{10}(Spin(9))\cong\Z_8$ 
\cite[p. 132, Table]{Mi}, we obtain
$$
\Delta_\K(i_\K\alpha_1(8))=b\alpha_1(7)=0
$$
and  $G_{11}(\mathbb{K}P^2)\supseteq\{i_\mathbb{K}\alpha_1(8)\}\cong\Z_3$. This leads to the equality $G_{11}(\mathbb{K}P^2)=P_{11}(\mathbb{K}P^2)=\{i_\mathbb{K}\alpha_1(8)\}$ and (2) follows.

We recall
$$
\pi_{15}(\K P^2)=\{i_K(E\sigma')\}\cong\Z_{120}.
$$
By 
(\ref{iKb}) and that $\pi_{14}(Spin(9))\cong\Z_8\oplus\Z_2$
\cite[p. 132, Table]{Mi}, we obtain
$$
\Delta_\K(i_\K\alpha_2(8))=0.
$$
Hence, we obtain\ $G_{15}(\mathbb{K}P^2)\ni i_\K\alpha_2(8)$ and
this leads to (3).

We recall \cite[Theorem 7.1, (7.4)]{T} $\pi_{15}(\S^7)=\{\sigma'\eta_{14},\bar{\nu}_7,\varepsilon_7\}$, where $\eta_7\sigma_8=\sigma'\eta_{14}+\bar{\nu}_7+\varepsilon_7$. By Lemma \ref{JT},
$$
\pi_{16}(\K P^2)={i_\K}_* E\{\sigma'\eta_{14},\bar{\nu}_7,\varepsilon_7\}.
$$
By (\ref{Kcyc}),
$i_\K(E\sigma')\eta_{15}\in P_{16}(\K P^2)$.
Next, by Lemma \ref{JT},
$[i_\K\bar{\nu}_8,i_\K]=i_\K(E\sigma')\bar{\nu}_{15}\ne 0$,  $[i_\K\varepsilon_8,i_\K]=i_\K(E\sigma')\varepsilon_{15}ne 0$ and
$[i_\K\bar{\nu}_8+\varepsilon_8,i_\K]=0$. Hence, by the relation $\eta_7\sigma_8=\sigma'\eta_7+\bar{\nu}_7+\varepsilon_7$, we have (4).

We recall
$$ \pi_{16}(\S^7)=\{\sigma'\eta^2_{14},\nu^3_7,\mu_7,\eta_7\varepsilon_8\}\cong(\Z_2)^4.
$$
By (4) and Lemma \ref{SV}, $i_\K(E(\sigma'\eta^2_{14})$ is cyclic.
By Lemma \ref{JT},
$$
[i_\K\nu^3_8,i_\K]=i_\K(E\sigma')\nu^3_{15}=i_\K\eta_8\bar{\varepsilon}_9
\ne 0, \ [i_\K\mu_8,i_\K]=i_\K(E\sigma')\mu_{15}\ne 0,\ [i_\K\nu^3_8,i_\K]=i_\K(E\sigma')
$$

$,\nu^3_7,\mu_7,\eta_7\varepsilon_8\}\cong(\Z_2)^4$.

This and Lemma \ref{JT} imply
$$
P_{17}(\K P^2)={i_\K}_* E\{\sigma'\eta^2_{14},\nu^3_7,\mu_7,\eta_7\varepsilon_8\}.
$$
We have $[\mu_8,\iota_8]=(E\sigma')\mu_{15}$. By Lemma \ref{JT},
$i_\K(E\sigma')\mu_{15}\ne 0$. So,
$i_\K\mu_8\not\in P_{17}(\K P^2)$. Moreover, by the relation $\sigma'\nu^3_{14}=\nu_7\sigma_{10}\nu^2_{17}=\eta_7\bar{\varepsilon}_8$ \cite[Lemma 12.10]{T}, we obtain
$i_\K(E\sigma')(\eta_{15}\varepsilon_{16})=i_\K(E\sigma')(\nu^3_{15})\ne 0$ and $i_\K(E\sigma')(\nu^3_{15}+\eta_{15})=0$ and this leads to (5).

We recall $\pi_{17}(\S^7)=\{\nu_7\sigma_{10},\eta_7\mu_8,\beta_1(7)\}$.
By Lemma \ref{JT},
$$
\pi_{18}(\K P^2)={i_\K}_*E\{\nu_7\sigma_{10},\eta_7\mu_8,\beta_1(7)\}\cong\Z_{24}\oplus\Z_2.
$$
Since $\pi_{17}(Spin(9))\cong\Z_8\oplus(\Z_2)^2$ \cite[p. 132, Table]{Mi},
$\Delta_\K(i_\K\beta_1(8))=b\beta_1(7)=0$ and
$i_\K\beta_1(8)\in G_{18}(\mathbb{K}P^2)$.
By (\ref{Kcyc}), $i_\K\nu_8\sigma_{11}$ is cyclic.

By Lemma \ref{JT} and the fact that $\pi_{24}(\mathbb{S}^7)=\{\sigma'\eta_{14}\mu_{15},\nu_7\kappa_{10},\bar{\mu}_7,\eta_7\mu_8\sigma_{17}\}\cong(\Z_2)^4$, we see that $[i_\K\eta_8\mu_9,i_\K]=i_\K(E\sigma')\eta_{15}\mu_{16}\neq 0$ in $\pi_{25}(\mathbb{K}P^2)$. Hence, $i_\K\eta_8\mu_9\not\in P_{18}(\K P^2)$
and this leads to (6).

Since $\pi_{20}(\K P^2)\cong\pi_{19}(\S^7)=0$ and $\pi_{19}(Spin(9))\cong\pi_{19}(F_4)\cong\Z_2$ \cite[p. 132, Table]{Mi}, we see that $\Delta_\K: \pi_{19}(\K P^2)\to\pi_{18}(Spin(9))$ is a monomorphism.
We recall
$$
\pi_{18}(\S^7)=\{\zeta_7,\bar{\nu}_7\nu_{10}\}\cong\Z_{504}\oplus\Z_2.
$$
We have $[\iota_8,\zeta_8]=(E\sigma')\zeta_{15}=x\zeta_8\sigma_{19}$ for $x$ odd \cite[Lemma 12.12]{T}. By Lemma \ref{JT} for $n=27$, we see that $i_\K(\zeta_8\sigma_{19})$ is of order $8$.
Since $\bar{\nu}_9\nu_{17}=[\iota_9,\nu_9]$ \cite[(7.22)]{T}, $i_\K\bar{\nu}_8\nu_{16}$ is cyclic by (\ref{Kcyc})
and this leads to (7).

By Lemma \ref{JT} and the fact that $\pi_{20}(\mathbb{S}^7)=\{\nu_7\sigma_{10}\nu_{17},\alpha_1(7)\beta_1(10)\}\cong\Z_6$, we obtain $\pi_{21}(\K P^2)={i_\K}_\ast\{\nu_8\sigma_{11}\nu_{18},\alpha_1(8)\beta_1(11)\}\cong\Z_6.$
By (\ref{iKb}),
$$
\Delta_\K(i_\K\alpha_1(8)\beta_1(11))=\Delta_\K(i_\K\alpha_1(8))\beta_1(10)=0.
$$
Hence, 
we obtain $G_{21}(\mathbb{K}P^2)\supseteq 2\pi_{21}(\K P^2)$.
By (6) and Lemma \ref{SV}, $i_\K(\nu_8\sigma_{11}\nu_{18})$ is cyclic.
This leads to (8) and completes the proof.
\end{pf}

At the end, write
$$
G''_k(\mathbb{F}P^n)=G_k(\mathbb{F}P^n)\cap({i_\mathbb{F}}_\ast E\pi_{k-1}(\S^{d-1})).
$$
Notice that $G''_k(\mathbb{F}P^n)=0$ if $d=1,2$ and $k\ge d+1$.
We close the paper with:

\begin{quest}
{\em What about $G''_k(\mathbb{H}P^n)$?}
\end{quest}

{\bf Acknowledgments.} The authors would like to express their thanks to Department of Mathematics,
Dalhousie University, Halifax (Canada) for its hospitality and support on June 29 -  July 05, 2008.
\par Further, the second author would like to thank for hospitalities and supports to Departments of Mathematics
and other sciences of Nicolaus Copernicus University, Toru\'n (Poland) on August 18 - 31, 2005
and  Korea University, Seoul (Korea) on September 25 - October 01, 2009.

\end{document}